\batchmode
\documentclass[11pt]{article}

\usepackage{epsfig}
\usepackage{graphicx}
\usepackage{color}
\usepackage{mathtools}

\newtheorem{theorem}{Theorem}
\newtheorem{lemma}{Lemma}
\newtheorem{corollary}{Corollary}

\newcommand{\be}{\begin{equation}}
\newcommand{\ee}{\end{equation}}
\newcommand{\bea}{\begin{eqnarray}}
\newcommand{\eea}{\end{eqnarray}}
\newcommand{\beas}{\begin{eqnarray*}}
\newcommand{\eeas}{\end{eqnarray*}}
\newcommand{\ba}{\begin{array}}
\newcommand{\ea}{\end{array}}

\DeclarePairedDelimiter{\floor}{\lfloor}{\rfloor}
\DeclarePairedDelimiter{\ceil}{\lceil}{\rceil}

\definecolor{armygreen}{rgb}{0.29, 0.33, 0.13}

\newcommand{\real}{\mbox{$\mathbb{R}$}}
\newcommand{\Natural}{\mbox{$\mathrm{I\!N}$}}

\newcommand{\eps}{\ensuremath{\epsilon}}

\newcommand{\bfp}{\ensuremath{\mathbf{p}}}

\newcommand{\bfx}{\ensuremath{\mathbf{x}}}
\newcommand{\bfy}{\ensuremath{\mathbf{y}}}

\newcommand{\bfF}{\ensuremath{\mathbf{F}}}

\newcommand{\bfX}{\ensuremath{\mathbf{X}}}

\def\XXint#1#2#3{{\setbox0=\hbox{$#1{#2#3}{\int}$}
     \vcenter{\hbox{$#2#3$}}\kern-.5\wd0}}

\newcommand{\alpot}{\ensuremath{\frac{\alpha}{2}}}

\newcommand{\mcF}{\ensuremath{\mathcal{F}}}
\newcommand{\mcG}{\ensuremath{\mathcal{G}}}
\newcommand{\mcH}{\ensuremath{\mathcal{H}}}

\newcommand{\mcV}{\ensuremath{\mathcal{V}}}

\newcommand{\wtilde}[1]{\ensuremath{\widetilde{#1}}}

\def\qed{\hbox{\vrule width 6pt height 6pt depth 0pt}}

\usepackage{amsmath}
\usepackage{amssymb}

\title{Characterization of the weighted Sobolev space $H_{\beta}^{s}(\Omega)$ in $\real^{2}$ 
in terms of the decay rate of Fourier-Jacobi coefficients} 
\author{	
	V.J.~Ervin\thanks{School of Mathematical and Statistical Sciences,
	  Clemson University, Clemson, South Carolina 29634-0975, USA.
	  email: {\tt vjervin@clemson.edu}. } }

\date{\today}

\begin{document}
\maketitle

\begin{abstract}
In this paper, motivated by the analysis of the fractional Laplace equation on the unit disk
in $\real^{2}$, we establish a characterization of the weighted Sobolev space $H_{\beta}^{s}(\Omega)$
in terms of the decay rate of Fourier-Jacobi coefficients. This framework is then used to give a precise
analysis of the solution to the fractional Laplace equation on the unit disk.

\end{abstract}



\setcounter{equation}{0}
\setcounter{figure}{0}
\setcounter{table}{0}
\setcounter{theorem}{0}
\setcounter{lemma}{0}
\setcounter{corollary}{0}
\setcounter{definition}{0}

\section{Introduction}
 \label{sec_intro}
Weighted Sobolev spaces play a fundamental role in the analysis of problems involving functions
with known singularities. The weighted spaces, by exploiting the known singularities, enable a precise
characterization of the functions involved and the mapping properties of the operators. Often, as
in the case studied herein, the functions of interest exhibit a specific singular behavior at the boundary
of the domain.

The investigation presented in this paper is motivated by the solution of the fractional Laplace equation
on the unit disk in $\real^{2}$. The solution is well known to exhibit a singular behavior at the boundary 
of the domain, where the normal derivative of the solution becomes unbounded \cite{aco171, ros141}

In \cite{dyd171}, Dyda, Kuznetsov, and Kwa\'{s}nicki showed that for $\omega(\bfx) \ := \ (1 \, - \| \bfx \|^{2})_{+}$,
$\mcV_{l, \mu}(\bfx)$ a solid harmonic polynomial and $P_{n}^{(\alpot , l)}(2 | \bfx |^{2} \, - \, 1)$ a Jacobi polynomial,
that the fractional Laplacian has the pseudo eigenfunction  (for $| \bfx | < 1$)
\[
  ( - \Delta )^{\alpot} \omega(\bfx)^{\alpot} \, \mcV_{l, \mu}(\bfx) \, P_{n}^{(\alpot , l)}(2 | \bfx |^{2} \, - \, 1) \ = \
  \lambda_{l, n} \, \mcV_{l, \mu}(\bfx) \, P_{n}^{(\alpot , l)}(2 | \bfx |^{2} \, - \, 1) \, , 
\]
for $\lambda_{l, n}$ a constant (see \eqref{bceq0v}).  This property was used by Hao, Li, Zhang, and Zhang in
\cite{hao211} to investigate the solution of the fractional Laplace equation on the unit disk within the framework of
a two parameter family of function spaces $\mathbf{B}_{\gamma}^{s_{1} , s_{2}}(\Omega)$, defined in terms of the
decay rate of Fourier-Jacobi coefficients. They were able to establish that for a suitable RHS function 
$f \in \mathbf{B}^{s_{1} , s_{2}}_{\alpot}(\Omega) \cap H^{-\alpot}(\Omega)$ a unique solution 
$\wtilde{u}(\bfx) \ = \ \omega(\bfx)^{\alpot} \, u(\bfx)$ exists, with 
$u \in \mathbf{B}^{s_{1} + \alpot \, , \, s_{2} + \alpha}_{\alpot}(\Omega)$.
In this paper, by identifying a more appropriate function setting that $\mathbf{B}_{\gamma}^{s_{1} , s_{2}}(\Omega)$,
we are able to 
give a more precise solution statement (see Corollary \ref{corext1}).

In \cite{bab011} Babu\v{s}ka and Guo established an equivalence between a family of weighted Sobolev spaces
(describing functions with singularities at the endpoints of an interval) with a family of function spaces
defined by the decay rate of a function's Fourier-Jacobi coefficients. This function framework was used 
in \cite{aco181, erv191, hao201, hao202, jia181, mao161, mao181, zhe211, zhe231} 
to study the existence and uniqueness of the solution to the fractional Laplace equation
on a bounded interval in $\real^{1}$, and its spectral approximation.

In this paper we first extend the result of Babu\v{s}ka and Guo \cite{bab011}, on the equivalence 
between a family of weighted Sobolev spaces and a family of function spaces
defined by the decay rate of a function's Fourier-Jacobi coefficients,
from the interval $(-1 , 1)$ in $\real^{1}$, to the
unit disk in $\real^{2}$. Then, using the framework provided by the family of functions 
defined by the decay rate of their Fourier-Jacobi coefficients, we provide a precise analysis of the existence and
uniqueness of the fractional Laplace equation on the unit disk.

Following in Sections \ref{sec_prelim}  and \ref{sec_FunS}  we introduce definitions and notation, 
and cite some important results 
used throughout the paper. In Sections \ref{ssec_eqH1}  and \ref{ssec_eqHs} we establish the equivalent between the
families of function spaces. Using the framework of the function spaces defined by the decay rate of the 
Fourier-Jacobi coefficients, in Section \ref{ssec_exHsneg} we give a precise characterization of the dual spaces.
Section \ref{sec_ExUng} presents the analysis of the fractional Laplace equation on the unit disk.
 

 \setcounter{equation}{0}
\setcounter{figure}{0}
\setcounter{table}{0}
\setcounter{theorem}{0}
\setcounter{lemma}{0}
\setcounter{corollary}{0}
\setcounter{definition}{0}
\section{Preliminaries}
 \label{sec_prelim}

In this section we present definitions and notation, and some results used later in the paper. In particular, we introduce
an orthogonal basis for the weighted $L^{2}$ space, $L^{2}_{\beta}(\Omega)$.

\subsection{Jacobi Polynomials}
\label{ssec_jpoly}
The Jacobi polynomials are defined in terms of the regularized Hypergeometric function $_{2}\bfF_{1}(\cdot)$ 
(see \cite{dyd171}) as 
\begin{align*}
P_{n}^{(a , b)}(t) &:= \ \frac{\Gamma(a \, + \, 1 \, + \, n)}{n !} \,
 _{2}\bfF_{1}\left( \begin{array}{c}
                          -n , \ 1 + a + b + n \\  a + 1  \end{array} \, \Big{\vert} \, \frac{1 - t}{2} \right)   \\ 
&= \  \frac{(-1)^{n} \, \Gamma(b \, + \, 1 \, + \, n)}{n !} \,
 _{2}\bfF_{1}\left( \begin{array}{c}
                          -n , \ 1 + a + b + n \\  b + 1  \end{array} \, \Big{\vert} \, \frac{1 + t}{2} \right)  \\  
&= \ (-1)^{n} \frac{\Gamma(n + 1 + b)}{n ! \, \Gamma(n + 1 + a + b)} \, \sum_{j = 0}^{n} (-1)^{j} \, 2^{-j} \, 
\left( \begin{array}{c}
           n  \\ j  \end{array} \right) \frac{\Gamma(n + j + 1 + a + b)}{\Gamma(j + 1 + b)} \, (1 + t)^{j} \, .  
\end{align*}


In case $a, \, b > -1$ the Jacobi polynomials satisfy the following orthogonality property.
\begin{align}
 & \int_{-1}^{1} (1 - t)^{a} (1 + t)^{b} \, P_{j}^{(a , b)}(t) \, P_{k}^{(a , b)}(t)  \, dt 
 \ = \
   \left\{ \begin{array}{ll} 
   0 , & k \ne j  \\
   |\| P_{j}^{(a , b)} |\|^{2}
   \, , & k = j  
    \end{array} \right.  \, ,  \nonumber \\
& \quad \quad \mbox{where } \  \ |\| P_{j}^{(a , b)} |\| \ = \
 \left( \frac{2^{(a + b + 1)}}{(2j \, + \, a \, + \, b \, + 1)} 
   \frac{\Gamma(j + a + 1) \, \Gamma(j + b + 1)}{\Gamma(j + 1) \, \Gamma(j + a + b + 1)}
   \right)^{1/2} \, .
  \label{spm22}
\end{align}

\subsection{Solid Harmonic Polynomials}
\label{ssec_shpoly}
The solid harmonic polynomials in $\real^{d}$ are the polynomials in $d$ variables which satisfy 
Laplace's equation.

In $\real^{2}$ the solid harmonic polynomials of degree $l$ can be conveniently written in polar
coordinates
as $\{ r^{l} \cos (l \varphi ) \ , \   r^{l} \sin (l \varphi ) \} $.

\subsection{The Function Space $L^{2}_{\beta}(\Omega)$}
\label{ssec_funspc}     
With respect to the Cartesian coordinate system, we denote a point $\bfx \in \real^{2}$ as $\bfx = (x , y)$, and
with respect to the polar coordinate system, $\bfx = (r , \varphi)$.

We let $\Omega$ denote the unit disk in $\real^{2}$, i.e., $\Omega \, := \, \{ \bfx = (x , y) \, : \, x^{2} + y^{2} < 1\} . $

For the weight function $w(\cdot)$, $w(\bfx) > 0 \, ,
\bfx \in \Omega$, the associated weighted $L^{2}$ function space is defined as
\[
 L^{2}_{w}(\Omega) \, := \, \{ f \, : \, \| f \|_{L^{2}_{w}} < \infty \} \, , \ \ \mbox{ where }
  \| f \|_{L^{2}_{w}}^{2} \, := \, \int_{\Omega} \, w(\bfx) \, ( f(\bfx) )^{2} \, d\Omega \, .
\]
Associated with $ L^{2}_{\beta}(\Omega) $ we have the inner product, defined for $f, \, g \, : \, \Omega \rightarrow \real^{n}$, by
\[
  (f , g)_{L^{2}_{w}} \, := \, \int_{\Omega} w(\bfx) \, f(\bfx) \cdot g(\bfx) \, d\Omega \, .
\]

The usual $L^{2}(\Omega)$ norm and inner product (corresponding to $w(\bfx) = 1$) are respectively denoted by
$\| f \|$ and $(f , g)$.

For $\bfx = (r , \varphi) \in \Omega$ , the weight functions
$ \omega^{\gamma} \, := \, (1 - r^{2})^{\gamma} $, for $\gamma \in \real$, play a central role in the analysis.
To simplify notation, for lower case Greek letters, $\alpha, \beta, \gamma \in \real$ we let
\be
 (f , g)_{\beta} \, := (f , g)_{L^{2}_{\omega^{\beta}}} \, = \, \int_{\Omega} \omega^{\beta}  \, f(\bfx) \cdot g(\bfx) \, d\Omega
 \ = \ \int_{\Omega}(1 - r^{2})^{\beta}  \, f(\bfx) \cdot g(\bfx) \, d\Omega \, ,
\label{deflb}
\ee
\[  
 \| f \|_{\beta} \, := \, (f , f)_{\beta}^{1/2} \ \ \mbox{and } \ \ L^{2}_{\beta}(\Omega) \, := \, L^{2}_{\omega^{\beta}}(\Omega) \, .
\]

Also, for notation convenience, let $\rho := \, 2r^{2} - 1$. \\

A basis for 
$L^{2}_{\beta}(\Omega)$ is given by a product of the solid harmonic polynomials and Jacobi polynomials.

\subsubsection{Basis for $L^{2}_{\beta}(\Omega)$}
\label{sssec_R2}     
\[
\mbox{Let } \ \mcV_{0 , 1}(\bfx) := \ \frac{1}{2} \, , \ 
\mcV_{l , 1}(\bfx) \ := \ r^{l} \, \cos(l \varphi) \, , \  \ l = 1, 2, \ldots  \ \ \mbox{ and } \
\mcV_{l , -1}(\bfx) \ := \ r^{l} \, \sin(l \varphi) \, , \  \ l = 1, 2, \ldots  .
\]

We use the following notation
\[
    \mcV_{l , \mu^{*}}(\bfx) \ = \ \left\{ \begin{array}{rl}
    \mcV_{l , -1}(\bfx) & \mbox{  if } \mu = 1 \, ,   \\
    \mcV_{l , 1}(\bfx) & \mbox{  if } \mu = -1 \, .  \end{array} \right.
\]

Additionally, for a linear operator $\mcF( \cdot )$,  we use $\mcF(\mcV_{l , \mu}(\bfx)) \ = \ (\pm) \mcV_{l , \sigma}(\bfx)$ to denote
\[
\mcF(\mcV_{l , \mu}(\bfx)) = \ (\pm) \mcV_{l , \sigma}(\bfx) \ =  \ \left\{ \begin{array}{rl}
   + \mcV_{l , \sigma}(\bfx) & \mbox{  if } \mu = 1 \, , \\
   - \mcV_{l , \sigma}(\bfx) & \mbox{  if } \mu = -1 \, .  \end{array} \right.
\]
For example,
\[
 \frac{\partial}{\partial \varphi} \mcV_{l , \mu}(\bfx) \ = \ (\mp) \mcV_{l , \mu^{*}}(\bfx) \, .
\]

An orthogonal basis for $L^{2}_{\beta}(\Omega)$ is  \cite{li141, wun051} (recall $\rho \, = \, 2 r^{2} \, - \, 1$)
\be
 \left\{  \cup_{l = 0}^{\infty} \cup_{n = 0}^{\infty} \left\{ \mcV_{l , 1}(\bfx) \, P_{n}^{(\beta  ,  l)}(\rho) \right\} \right\}   \ \cup
  \left\{ \cup_{l = 1}^{\infty} \cup_{n = 0}^{\infty}  \left\{ \mcV_{l , -1}(\bfx) \, P_{n}^{(\beta  ,  l)}(\rho) \right\} \right\}    \, .
\label{bsR2}
\ee

For notation brevity  we denote the basis in \eqref{bsR2} as
\[
  \cup_{l = 0}^{\infty} \cup_{n = 0}^{\infty}  \cup_{\mu = {1 , -1}} \left\{ \mcV_{l , \mu}(\bfx) \, P_{n}^{(\beta  ,  l)}(\rho) \right\} 
\]
where we implicit assume that the terms   $\mcV_{0, -1}(\bfx) \, P_{n}^{(\beta  ,  l)}(\rho), \, n = 0, 1, \ldots$ are omitted 
from the set. 

We trivially extend the definitions of $\mcV_{l, \mu}(\bfx)$ and $P_{n}^{(\alpha  ,  \beta)}(t)$ to negative integer values for $l$ and $n$
by $\mcV_{l, \mu}(\bfx) = 0$ for $l \le -1$, and $P_{n}^{(\alpha  ,  \beta)}(t) = 0$ for $n \le -1$.

In \cite{zhe232} we used the notation $V_{l, \mu}(\bfx)$ to represent the harmonic polynomials, with $\mcV_{l, \mu}(\bfx) = V_{l, \mu}(\bfx)$,
except for $\mcV_{0, 1}(\bfx) = \, \frac{1}{2} V_{0, 1}(\bfx) = \frac{1}{2}$. This change of notation, together with the extension described in
the previous paragraph allows a more compact representation for the 
following derivative formulas.


\begin{corollary} \cite{zhe232} \label{gradwR2}
 Let $f(\bfx) \ = \ \mcV_{l , \mu}(\bfx) P_{n}^{(\gamma , l)}(2r^{2} - 1)$. Then,
 \begin{align}
\frac{\partial f}{\partial x} &= \
  (n + l)  \, \mcV_{l-1 , \mu}(\bfx) P_{n}^{(\gamma + 1 \, , \, l - 1)}(2r^{2} - 1) 
 \ + \
(n + \gamma + l + 1) \,  \mcV_{l+1 , \mu}(\bfx) P_{n-1}^{(\gamma + 1 \, , \, l + 1)}(2r^{2} - 1) \, ,    \label{dersw1} \\
\mbox{and} &  \nonumber \\
\frac{\partial f}{\partial y} &= \
 (\mp) (n + l) \,  \mcV_{l-1 \, , \mu^{*}}(\bfx) \, P_{n}^{(\gamma + 1 \, , \,  l - 1)}(2 r^{2} - 1)    \nonumber \\
 & \hspace{2.0in} 
\ + \   (\pm) (n + \gamma + l + 1) \, \mcV_{l+1 \, , \mu^{*}}(\bfx) \,   P_{n-1}^{(\gamma + 1 \,  ,  \, l+1)}(2 r^{2} - 1)  \, .  
  \label{dersw2}  
 \end{align}
 \end{corollary}
    
  
For notation convenience,  we let $\mathbb{N}_{0}  := \mathbb{N} \cup {0}$ and
 use $a \sim b$ to denote that there exists constants $C_{0}$ and $C_{1} > 0$ such that 
 $C_{0} \, b \, \le \, a \,  \le \, C_{1} \, b$. Additionally, we use $a \, \lesssim \, b$ to denote that there exists a constant $C_{1} > 0$ such that
 $a \, \le \, C_{1} \,  b$. 
 
 For $s \in \mathbb{R}$, $\floor{s}$ is used to denote the largest integer that is less than or equal to $s$, and
 $\ceil{s}$ is used to denote the smallest integer that is greater than or equal to $s$.
 
From Stirling's formula we have that
\begin{equation}
\lim_{n \rightarrow \infty} \, \frac{\Gamma(n + \sigma)}{\Gamma(n) \, n^{\sigma}}
\ = \ 1 \, , \mbox{ for } \sigma \in \mathbb{R}.  
 \label{eqStrf}
\end{equation} 
 

 \setcounter{equation}{0}
\setcounter{figure}{0}
\setcounter{table}{0}
\setcounter{theorem}{0}
\setcounter{lemma}{0}
\setcounter{corollary}{0}
\setcounter{definition}{0}
\section{The function spaces $H_{\gamma}^{s}(\Omega)$ and $\wtilde{H}_{\gamma}^{s}(\Omega)$}
\label{sec_FunS}

In this section we introduce the function spaces used in the analysis below.
The first space defined, $H_{\gamma}^{s}(\Omega)$, is based upon the boundedness
of weighted $L^{2}$ norms of a function derivatives. 
The second space defined, $\wtilde{H}_{\gamma}^{s}(\Omega)$, is based upon the decay 
rates of the Fourier-Jacobi coefficients of a function in $L^{2}_{\gamma}(\Omega)$.

\subsection{The space $H_{\gamma}^{s}(\Omega)$}
\label{ssec_H01}

Similar to Bab\v{u}ska and Guo \cite{bab011} and Guo and Wang \cite{guo041}, define the 
weighted Sobolev space $H_{\gamma}^{s}(\Omega)$, for $s \in \Natural$, as
\be
H_{\gamma}^{s}(\Omega) \, := \, \Big\{ f \in L^{2}_{\gamma}(\Omega) \, : \, 
| f |_{H_{\gamma}^{s}} \, := \, \left( \sum_{j = 0}^{s} \left( \begin{array}{c}   s  \\ j  \end{array} \right) \, 
\left\| \frac{\partial^{s} f(\bfx)}{\partial y^{j} \, \partial x^{s - j}} \right\|^{2}_{\gamma + s} \right)^{1/2} \, < \, \infty \Big\} \, .
\label{defHsp}
\ee
We associate with $H_{\gamma}^{s}(\Omega)$ the norm
\[
  \| f \|_{H_{\gamma}^{s}} \, := \, \left( \| f \|_{\gamma}^{2} \ + \  | f |_{H_{\gamma}^{s}}^{2} \right)^{1/2} \, .
\]

For $s > 0$, $s \not\in \Natural$,  $H_{\gamma}^{s}(\Omega)$ is defined using the K-method of interpolation. For
$s < 0$ the spaces are defined by (weighted) $L^{2}_{\gamma}$ duality.

These function spaces are particularly useful when dealing with functions that have a known singular behavior 
at the boundary of the domain, such as the solution to fractional Laplace equations.

\subsection{The space $\wtilde{H}_{\gamma}^{s}(\Omega)$ }
\label{ssec_H02}

Assume $\gamma \geq 0$ is given. (Recall that $\rho \, = \, 2r^{2} - 1$.) Let
\be
u(\bfx) \ = \ \sum_{l, n, \mu} a_{l, n, \mu} \, \mcV_{l, \mu}(\bfx) \, P_{n}^{(\gamma, l)}(\rho) \, ,  \ \  \mbox{and  }
h_{l, n} \ := \ \| \mcV_{l, \mu}(\bfx) \, P_{n}^{(\gamma, l)}(\rho) \|_{\gamma} \, .
\label{defunh}
\ee

For $s \geq 0$, the space $\wtilde{H}_{\gamma}^{s}(\Omega)$ is defined as
\be
\wtilde{H}_{\gamma}^{s}(\Omega) \ := \ \Big\{ u(\cdot) \in L^{2}_{\gamma}(\Omega) \, : \, 
\sum_{l, n, \mu} (n + 1)^{s} (n + l + 1)^{s} \, a_{l, n, \mu}^{2} \, h_{l, n}^{2}  < \infty \Big\} \, ,
\label{defXs}
\ee
with 
\be
\| u \|_{\wtilde{H}_{\gamma}^{s}}  := \ \Big( \sum_{l, n, \mu} (n + 1)^{s} (n + l + 1)^{s} \, a_{l, n, \mu}^{2} \, h_{l, n}^{2} \Big)^{1/2}.
\label{nmHt}
\ee

These function spaces are particularly useful when dealing with functions expressed as a linear combination
of basis functions of $L^{2}_{\beta}(\Omega)$.

A first step toward establishing the equivalence of $H_{\gamma}^{s}(\Omega)$ and $\wtilde{H}_{\gamma}^{s}(\Omega)$
is to show that
the family of spaces $\wtilde{H}_{\gamma}^{s}(\Omega)$ is a family of interpolation spaces.

\begin{theorem} \label{thmXsip}
The spaces $\wtilde{H}_{\gamma}^{s}(\Omega)$ defined in \eqref{defXs} is a family of interpolation spaces with 
respect to the $K-$method of interpolation.
\end{theorem}

\textbf{Outline of the proof}: To establish that $\wtilde{H}_{\gamma}^{s}(\Omega)$  is a family of interpolation spaces with 
respect to the $K-$method of interpolation, we need to show that for 
for nonnegative integers $k$ and $m$, with $k < s < m$, specifically $s \, = \, (1 - \theta) k \, + \, \theta \, m$ 
for some $0 < \theta < 1$,
\be
\| u \|_{\wtilde{H}_{\gamma}^{s}} \ = \ \left( \int_{0}^{\infty} t^{-2 \theta} \, | K(t, u) |^{2} \, \frac{dt}{t} \right)^{1/2} \, , \ \ 
\mbox{where } \ K(t, u) \ = \ 
\inf_{v \in \wtilde{H}_{\gamma}^{m}} \left( \|u - v \|_{\wtilde{H}_{\gamma}^{k}} \ 
+  \ t \, \| v \|_{\wtilde{H}_{\gamma}^{m}} \right) \, .
\label{intnrm1}
\ee
Noting that 
\[
K(t, u) \approx \wtilde{K}(t, u) \ = \ \inf_{v \in \wtilde{H}_{\gamma}^{m}} 
\left( \|u - v \|_{\wtilde{H}_{\gamma}^{k}}^{2} \ +  \ t^{2} \, \| v \|_{\wtilde{H}_{\gamma}^{m}}^{2} \right)^{1/2} \, ,
\]
for $u$ given by \eqref{defunh},
we determine $\wtilde{K}(t, u)$ explicitly as a function of $t$, $a_{l, n, \mu}$, and $h_{l, n}$.
We then show that $\| u \|_{\wtilde{H}_{\gamma}^{s}}$ given by \eqref{intnrm1} is equivalent to 
$\| u \|_{\wtilde{H}_{\gamma}^{s}}$ given by \eqref{nmHt}.

The lengthy proof
is presented in Appendix \ref{sec_Intprf}.  \\
\mbox{  } \hfill \qed

 \setcounter{equation}{0}
\setcounter{figure}{0}
\setcounter{table}{0}
\setcounter{theorem}{0}
\setcounter{lemma}{0}
\setcounter{corollary}{0}
\setcounter{definition}{0}
\section{Equivalence of the spaces $H_{\gamma}^{1}(\Omega)$ and $\wtilde{H}_{\gamma}^{1}(\Omega)$}
\label{ssec_eqH1}

To indicate how the norms of the spaces $H_{\gamma}^{s}(\Omega)$ and $\wtilde{H}_{\gamma}^{s}(\Omega)$
are equivalent, we begin by showing the norm 
equivalence for $H_{\gamma}^{1}(\Omega)$ and $\wtilde{H}_{\gamma}^{1}(\Omega)$.

Using Corollary \ref{gradwR2} it follows fairly easily (using orthogonality of the basis functions and the triangle inequality) that
$ \| u \|_{H_{\gamma}^{1}} \lesssim \| u \|_{\wtilde{H}_{\gamma}^{1}}$. The challenging part is to establish that 
$ \| u \|_{H_{\gamma}^{1}} \gtrsim \| u \|_{\wtilde{H}_{\gamma}^{1}}$. This difficulty
is due to the occurrence of two different basis functions
on the RHS of \eqref{dersw1} and \eqref{dersw2}. However a surprising cancellation occurs when 
$\| \frac{\partial u}{\partial x} \|_{L^{2}_{\gamma + 1}}^{2}$ and $\| \frac{\partial u}{\partial y} \|_{L^{2}_{\gamma + 1}}^{2}$ 
are added together (see \eqref{poit3p5} and \eqref{poit5p5}) below). This cancellation is similar to that which occurs in 
the addition: $(a_{1} \, + \, a_{2})^{2} \ + \ (a_{1} \, - \, a_{2})^{2} \ = \ 2 ( a_{1}^{2} \, + \, a_{2}^{2} )$.

We first establish an expression for $| \cdot |_{H_{\gamma}^{1}} $.

\begin{lemma}  \label{snH01}
For $u \in H_{\gamma}^{1}(\Omega)$,
\be
| u |_{H^{1}_{\gamma}}^{2} \ \sim \  \sum_{2n + l \ge 1, \, \mu} (n + 1) \, (n + l + 1) \, a_{l, n, \mu}^{2} \, 
  \| \mcV_{l, \mu}(\bfx) P_{n}^{(\gamma , l)}(\rho) \|_{\gamma}^{2} \, . 
\label{repsn1}
\ee
\end{lemma}
In terms of notation simplicity, there are two issues which arise in the proof. The equality of 
$ \| \mcV_{l , \mu}(\bfx) P_{n}^{(\gamma + 1 \, , \, l)}(\rho) \|_{\gamma + 1} \ = \ 
 \| \mcV_{l , \mu^{*}}(\bfx) P_{n}^{(\gamma + 1 \, , \, l)}(\rho) \|_{\gamma + 1} $ for 
 $(l , \, \mu)~\neq~(0 , \, \pm 1)$ is used in \eqref{poit3p5} to combine terms in 
 $\left\| \frac{\partial u}{\partial x} \right\|^{2}_{\gamma + 1}$ and
 $\left\| \frac{\partial u}{\partial y} \right\|^{2}_{\gamma + 1} $. As $\mcV_{0, -1}(\bfx) \, = \, 0$ (formally not 
 a solid harmonic polynomial)  $\| \mcV_{0, -1}(\bfx) \, P_{n}^{(\gamma + 1, \, 0)}(\rho) \|_{\gamma + 1} \, \neq \, 
\| \mcV_{0, 1}(\bfx) \, P_{n}^{(\gamma + 1, \, 0)}(\rho) \|_{\gamma + 1}$, hence the expression in \eqref{poit3p5}
must be written as two summations to account for this exception. From Lemma \ref{lmaker3}, as
$\big\{ \mcV_{0, \mu} \, P_{0}^{(\gamma, \, 0)}(\rho) \big\} \, = \, Null(\mcG_{1})$, the coefficients $a_{0, 0, 1}$ and
$a_{0, 0, -1}$ do not formally appear in the expression for $| u |_{H^{1}_{\gamma}}^{2}$. However, for
notation simplicity we leave the terms in the analysis (with the knowledge that they can assumed to be zero) until
the final step. 

\textbf{Proof}: 
Let $u(\bfx) \ = \ \sum_{l, n, \mu} a_{l, n, \mu} \, \mcV_{l, \mu}(\bfx) \, P_{n}^{(\gamma , l )}(\rho)$.
From Corollary \ref{gradwR2},
\begin{align}
\frac{\partial u}{\partial x} &= \
\sum_{l, n, \mu} a_{l, n, \mu} \left( 
 (n + l)  \, \mcV_{l-1 , \mu}(\bfx) P_{n}^{(\gamma + 1 \, , \, l - 1)}(\rho) 
 \ + \
(n + \gamma + l + 1) \,  \mcV_{l+1 , \mu}(\bfx) P_{n-1}^{(\gamma + 1 \, , \, l + 1)}(\rho) \right)  \nonumber \\
&= \ \sum_{l, n, \mu} \underbrace{\left( a_{l+1 , n, \mu} (n + l + 1) \, + \, a_{l-1 , n+1 , \mu} (n + \gamma + l + 1) \right)}_{:= 
\, _{x}A_{l, n, \mu}}
 \, \mcV_{l , \mu}(\bfx) P_{n}^{(\gamma + 1 \, , \, l)}(\rho) \, ,   \label{poit2}  \\
&= \ \sum_{l, n, \mu} \,  _{x}A_{l, n, \mu}
 \, \mcV_{l , \mu}(\bfx) P_{n}^{(\gamma + 1 \, , \, l)}(\rho) \, .  \label{poit2p5}  
\end{align}

Similarly,
\begin{align}
\frac{\partial u}{\partial y} &= \
 \ \sum_{l, n, \mu} \underbrace{\left( (\mp) a_{l+1 , n, \mu} (n + l + 1) \, + \, (\pm) a_{l-1 , n+1 , \mu} (n + \gamma + l + 1) \right)}_{:= 
\, _{y}A_{l, n, \mu}}
 \, \mcV_{l , \mu^{*}}(\bfx) P_{n}^{(\gamma + 1 \, , \, l)}(\rho) \, ,   \label{poit4}  \\
&= \ \sum_{l, n, \mu} \, _{y}A_{l, n, \mu} \, \mcV_{l , \mu^{*}}(\bfx) P_{n}^{(\gamma + 1 \, , \, l)}(\rho) \, . \label{poit5}
\end{align}

Hence, using orthogonality, $\mcV_{0 , -1}(\bfx) = 0$,
 and that $ \| \mcV_{l , \mu}(\bfx) P_{n}^{(\gamma + 1 \, , \, l)}(\rho) \|_{\gamma + 1} \ = \ 
 \| \mcV_{l , \mu^{*}}(\bfx) P_{n}^{(\gamma + 1 \, , \, l)}(\rho) \|_{\gamma + 1} $ for 
 $(l , \, \mu)~\neq~(0 , \, \pm 1)$
 
\begin{align}
& | u |_{H^{1}_{\gamma}(\Omega)}^{2} \ = \ 
\left\| \frac{\partial u}{\partial x} \right\|^{2}_{\gamma + 1} \ + \  \left\| \frac{\partial u}{\partial y} \right\|^{2}_{\gamma + 1} 
  \nonumber  \\
&= \ \sum_{n} \left( _{x}A_{0, n, 1}^{2}  \, + \, _{y}A_{0, n, -1}^{2}  \right)
 \, \| \mcV_{0 , 1}(\bfx) P_{n}^{(\gamma + 1 \, , \, 0)}(\rho) \|_{\gamma + 1}^{2} \, ,   \nonumber   \\
& \quad +  \sum_{l \ge 1, \, n, \mu} \left( _{x}A_{l, n, \mu}^{2}  \, + \, _{y}A_{l, n, \mu}^{2}  \right)
 \, \| \mcV_{l , 1}(\bfx) P_{n}^{(\gamma + 1 \, , \, l)}(\rho) \|_{\gamma + 1}^{2} \, ,   \label{poit3}   \\
&= \   \sum_{n}  \big( a_{1 , n, 1}^{2} \ + \ a_{1 , n, -1}^{2} \big) \, (n + 1)^{2}  
 \, \| \mcV_{0 , 1}(\bfx) P_{n}^{(\gamma + 1 \, , \, 0)}(\rho) \|^{2}_{\gamma + 1}   \nonumber \\
& \ \ + \sum_{l \ge 1, n, \mu} \Big( \big( a_{l+1 , n, \mu} (n + l + 1) \, + \, a_{l-1 , n+1 , \mu} (n + \gamma + l + 1) \big)^{2} 
 \nonumber \\
& \hspace{0.25in} + \, \big( (\mp) a_{l+1 , n, \mu} (n + l + 1) \, + \, (\pm) a_{l-1 , n+1 , \mu} (n + \gamma + l + 1) \big)^{2} \Big)
 \, \| \mcV_{l , 1}(\bfx) P_{n}^{(\gamma + 1 \, , \, l)}(\rho) \|^{2}_{\gamma + 1}   \label{poit3p5}  \\
&= \   \sum_{n}  \big( a_{1 , n, 1}^{2} \ + \ a_{1 , n, -1}^{2} \big) \, (n + 1)^{2}  
 \, \| \mcV_{0 , 1}(\bfx) P_{n}^{(\gamma + 1 \, , \, 0)}(\rho) \|^{2}_{\gamma + 1}   \nonumber \\
 & \ \ + \ 2 \,  \sum_{l \ge 1, \, n, \mu}  
 \Big( a_{l+1 , n, \mu}^{2} (n + l + 1)^{2} \, + \, a_{l-1 , n+1 , \mu}^{2} (n + \gamma + l + 1)^{2} \Big)
 \, \| \mcV_{l , 1}(\bfx) P_{n}^{(\gamma + 1 \, , \, l)}(\rho) \|^{2}_{\gamma + 1}   \, .   \label{poit5p5}
\end{align}

Using Lemma \ref{lmaeq3},
\be
 \| \mcV_{l , \mu}(\bfx) P_{n}^{(\gamma + 1 \, , \, l)}(\rho) \|^{2}_{L^{2}_{\gamma + 1}}  \ = \ 
  \frac{2n + \gamma + l + 1}{2n + \gamma + l + 2} \, \frac{n + \gamma + 1}{n + \gamma + l + 1} \,
  \| \mcV_{l , \mu}(\bfx) P_{n}^{(\gamma , l)}(\rho) \|^{2}_{L^{2}_{\gamma}}  \, .
 \label{poit6}
\ee
Also, from Lemma \ref{lmaeq4},
\be
 \| \mcV_{l , \mu}(\bfx) P_{n}^{(\gamma , l)}(\rho) \|^{2}_{L^{2}_{\gamma}}  \ = \  \frac{C_{l, \mu}}{C_{l+1 , \mu}} \, 
  \frac{2n + \gamma + l + 2}{2n + \gamma + l + 1} \, \frac{n + \gamma + l + 1}{n + l + 1} \,
  \| \mcV_{l+1 \,  , \, \mu}(\bfx) P_{n}^{(\gamma , \,  l + 1)}(\rho) \|^{2}_{L^{2}_{\gamma}}  \, .
 \label{poit7}
\ee
Additionally, from Lemmas \ref{lmaeq5},
\begin{align}
\| \mcV_{l , \mu}(\bfx) P_{n}^{(\gamma , l)}(\rho) \|^{2}_{L^{2}_{\gamma}}  
&= \     \frac{C_{l, \mu}}{C_{l-1 , \mu}} \, 
  \frac{2n + \gamma + l + 2}{2n + \gamma + l + 1} \, \frac{n + 1}{n + \gamma + 1} \,
  \| \mcV_{l-1 \,  , \, \mu}(\bfx) P_{n+1}^{(\gamma , \,  l - 1)}(\rho) \|^{2}_{L^{2}_{\gamma}}  \, .   \label{poit8}
\end{align}
 
Combining \eqref{poit3}, \eqref{poit5p5}--\eqref{poit8} we obtain
\begin{align}
& | u |_{H^{1}_{\gamma}}^{2} \ = \  
\sum_{l, n, \mu} \left( _{x}A_{l, n, \mu}^{2}  \, + \, _{y}A_{l, n, \mu}^{2}  \right)
 \, \| \mcV_{l , \mu}(\bfx) P_{n}^{(\gamma + 1 \, , \, l)}(\rho) \|_{\gamma + 1}^{2} \, ,   \label{poit8p3}   \\
& \ \ = \  
\sum_{n}  \big( a_{1 , n, 1}^{2} \ + \ a_{1 , n, -1}^{2} \big) \, \frac{C_{0, 1}}{C_{1, 1}} \, (n + 1) \, (n + \gamma + 1) \,   
 \, \| \mcV_{1 , 1}(\bfx) P_{n}^{(\gamma \, , \, 1)}(\rho) \|^{2}_{\gamma}   \nonumber \\
 & \quad + \sum_{l \ge 1, \, n, \mu} 2 \, \frac{C_{l, \mu}}{C_{l+1, \mu}} \, a_{l+1 , n, \mu}^{2} (n + \gamma + 1) \, (n + l + 1) \, 
  \| \mcV_{l+1 , 1}(\bfx) P_{n}^{(\gamma , \, l + 1)}(\rho) \|^{2}_{\gamma}    \nonumber  \\
& \hspace*{0.75in} + \   
 \sum_{l \ge 1, \, n, \mu} 2 \, \frac{C_{l, \mu}}{C_{l-1, \mu}} \, a_{l-1 , n+1, \mu}^{2} (n + 1) \, (n + \gamma + l + 1) \, 
  \| \mcV_{l-1 , 1}(\bfx) P_{n+1}^{(\gamma , \, l - 1)}(\rho) \|^{2}_{\gamma}     \label{poit8p4}   \\
& \ \ \sim \ \sum_{n}  \big( a_{1 , n, 1}^{2} \ + \ a_{1 , n, -1}^{2} \big) \, \frac{C_{0, 1}}{C_{1, 1}} \, (n + 1) \, (n + \gamma + 1) \,   
 \, \| \mcV_{1 , 1}(\bfx) P_{n}^{(\gamma \, , \, 1)}(\rho) \|^{2}_{\gamma}   \nonumber \\
 & \ +   \sum_{2n + l \ge 1, \, \mu}
\left( \frac{C_{l-1, \mu}}{C_{l, \mu}} \, (n + \gamma + 1) \, (n + l) \ + \ \frac{C_{l+1, \mu}}{C_{l, \mu}} \,
 n \, (n + \gamma + l + 1)  \right) \,
 a_{l, n, \mu}^{2} \, 
  \| \mcV_{l, 1}(\bfx) P_{n}^{(\gamma , l)}(\rho) \|^{2}_{\gamma}    \nonumber  \\
%
%
& \ \ \sim \  \sum_{2n + l \ge 1, \, \mu} (n + 1) \, (n + l + 1) \, a_{l, n, \mu}^{2} \, 
  \| \mcV_{l, \mu}(\bfx) P_{n}^{(\gamma , l)}(\rho) \|^{2}_{\gamma}  \, ,     \label{poit8p5} 
\end{align}
where in the next to the last step we have used Lemma \ref{lmaker3}.  Specifically, as 
$\partial / \partial x \big(  \mcV_{l, \mu}(\bfx) P_{n}^{(\gamma , l)}(\rho) \big) = 0$ and  
$\partial / \partial y \big(  \mcV_{l, \mu}(\bfx) P_{n}^{(\gamma , l)}(\rho) \big) = 0$ for $2n + l < 1$, the
coefficients $a_{l, n, \mu}$ vanish in \eqref{poit8p4} for $2n + l < 1$. 
Equivalently, we can assume that in  \eqref{poit8p4} $a_{l, n, \mu} = 0$ for $2n + l < 1$.  \\
\mbox{ } \hfill \qed

\begin{corollary}  \label{H1eqs}
The spaces $H^{1}_{\gamma}(\Omega)$ and $\wtilde{H}^{1}_{\gamma}(\Omega)$ are equivalent.
\end{corollary}
\textbf{Proof}: Given the definitions of the spaces $H^{1}_{\gamma}(\Omega)$ and $\wtilde{H}^{1}_{\gamma}(\Omega)$,
\eqref{defHsp} and \eqref{defXs}, we show that the spaces are equivalent by showing that their norms are equivalent.

For $u(\bfx) \ = \ \sum_{l, n, \mu} a_{l, n, \mu} \, \mcV_{l, \mu}(\bfx) \, P_{n}^{(\gamma , l )}(\rho) \in L^{2}_{\gamma}(\Omega)$, using 
orthogonality, we have
\be
 \| u \|_{\gamma}^{2} \ = \ \sum_{l, n, \mu} a_{l, n, \mu}^{2} \, \| \mcV_{l, \mu}(\bfx) \, P_{n}^{(\gamma , l )}(\rho) \|_{\gamma}^{2} \, .
 \label{poiy1}
\ee
Combining \eqref{poiy1} with \eqref{repsn1} yields
\begin{align*}
\| u \|_{H^{1}_{\gamma}}^{2} &\sim \ 
\sum_{l, n, \mu} a_{l, n, \mu}^{2} \, \| \mcV_{l, \mu}(\bfx) \, P_{n}^{(\gamma , l )}(\rho) \|_{\gamma}^{2} \ 
+  \sum_{2n + l \ge 1, \, \mu} (n + 1) \, (n + l + 1) \, a_{l, n, \mu}^{2} \, 
  \| \mcV_{l, \mu}(\bfx) P_{n}^{(\gamma , l)}(\rho) \|_{\gamma}^{2}  \\
&\sim \   \sum_{l, n, \mu} (n + 1) \, (n + l + 1) \, a_{l, n, \mu}^{2} \, 
  \| \mcV_{l, \mu}(\bfx) P_{n}^{(\gamma , l)}(\rho) \|_{\gamma}^{2}   \\
&\sim \   \| u \|_{\wtilde{H}^{1}_{\gamma}}^{2} \, .
\end{align*}
\mbox{  } \hfill \qed

 \setcounter{equation}{0}
\setcounter{figure}{0}
\setcounter{table}{0}
\setcounter{theorem}{0}
\setcounter{lemma}{0}
\setcounter{corollary}{0}
\setcounter{definition}{0}
\section{Equivalence of the spaces $H_{\gamma}^{s}(\Omega)$ and $\wtilde{H}_{\gamma}^{s}(\Omega)$}
\label{ssec_eqHs}

In this section we show that the spaces $H_{\gamma}^{s}(\Omega)$ and $\wtilde{H}_{\gamma}^{s}(\Omega)$ 
are equivalent. As, for $s \not\in \Natural$, both spaces are defined by the K-method of interpolation, all we need
to do is to establish the equivalence of the spaces for $s \in \Natural$. To show the equivalence of the space for
$s \in \Natural$ we just need to show that  the semi-norms 
$| \cdot |_{H_{\gamma}^{s}} $ and $| \cdot |_{\wtilde{H}_{\gamma}^{s}}  $ are equivalent.

The surprising cancellation noted at the beginning of Section \ref{ssec_eqH1} is again the key to show the equivalence 
of the $s$ semi-norms. Instead of considering each of the $s$ order derivative terms individually and then
summing to obtain the $s$ semi-norm, we consider the sum of the pair of terms 
$\| \frac{\partial p}{\partial x} \|_{L^{2}_{\gamma + s}}$ and $\| \frac{\partial p}{\partial y} \|_{L^{2}_{\gamma + s}}$
where $p$ is a $(s - 1)^{st}$ order derivative of $u$, and then sum over all the $p$. Doing so we are able
to utilize the noted cancellation property. In Appendix \ref{sec_snorm2} we prove the case for $s = 2$ before 
presenting the case for general $s \in \Natural$.

\begin{lemma}  \label{snH1}
For $u \in H_{\gamma}^{s}(\Omega)$, $s \in \Natural$, we have
\be
| u |_{H^{s}_{\gamma}}^{2} \ \sim \  \sum_{2n + l \ge s, \, \mu} (n + 1)^{s} \, (n + l + 1)^{s} \, a_{l, n, \mu}^{2} \, 
  \| \mcV_{l, \mu}(\bfx) P_{n}^{(\gamma , l)}(\rho) \|_{\gamma}^{2} \, . 
\label{repsns}
\ee
\end{lemma}
\textbf{Outline of the proof}: The proof is by induction.
As mentioned above, the key idea for the proof is to, instead of 
considering each of the $s$ order derivative terms individually and then
summing to obtain the $s$ semi-norm, we consider the sum of the pair of terms 
$\| \frac{\partial p}{\partial x} \|_{L^{2}_{\gamma + s}}$ and $\| \frac{\partial p}{\partial y} \|_{L^{2}_{\gamma + s}}$
where $p$ is a $(s - 1)^{st}$ order derivative of $u$, and then sum over all the $p$.
In Appendix \ref{sec_snorm2} we firstly show the result for $s = 2$ to illustrate the structure of the proof,
before discussing the general case for $s \in \Natural$. \\
\mbox{  } \hfill \qed

\begin{theorem}  \label{Hseqs}
For $s \ge 0$, the spaces $H^{s}_{\gamma}(\Omega)$ and $\wtilde{H}^{s}_{\gamma}(\Omega)$ are equivalent.
\end{theorem}
\textbf{Proof}:
For $s = 0$ the equivalence of the spaces follow immediately from their definitions. The case for $s = 1$ is established
in Corollary \ref{H1eqs}. For $s \in \Natural, \, s \ge 2$, the proof follows from \eqref{repsns} and the analogous argument
used for the proof of Corollary \ref{H1eqs}. For $s > 0, \, s \not\in \Natural$ the result follows from the equivalence of the spaces
for $s \in \Natural_{0}$ and the fact that the non integer order spaces are both defined using the K-method of interpolation. \\
\mbox{  } \hfill \qed

 \setcounter{equation}{0}
\setcounter{figure}{0}
\setcounter{table}{0}
\setcounter{theorem}{0}
\setcounter{lemma}{0}
\setcounter{corollary}{0}
\setcounter{definition}{0}
\section{Explicit characterization of the space $H_{\gamma}^{-s}(\Omega)$, for $s \in \real^{+}$}
\label{ssec_exHsneg}

In this section we obtain an explicit representation for the dual space of $H^{s}_{\gamma}(\Omega)$. Given 
the equivalence of the spaces established above, we define the inner product on $H^{s}_{\gamma}(\Omega)$
$\langle \cdot , \cdot \rangle_{H^{s}_{\gamma}}$, for 
$f(\bfx) \, = \, \sum_{l, n, \mu} f_{l, n, \mu} \, \mcV_{l, \mu}(\bfx) \, P_{n}^{(\gamma , l )}(\rho)$ and
$g(\bfx) \, = \, \sum_{l, n, \mu} g_{l, n, \mu} \, \mcV_{l, \mu}(\bfx) \, P_{n}^{(\gamma , l )}(\rho)$ as
\begin{align}
\langle f , g \rangle_{H^{s}_{\gamma}} &:= \,
\sum_{l, n, \mu} (n + 1)^{s} (n + l + 1)^{s} \, f_{l, n, \mu} \, g_{l, n, \mu} \, h_{l, n}^{2} \, , \   \nonumber \\ 
\mbox{ and use } \| f \|_{H^{s}_{\gamma}} &:= \langle f , f \rangle_{H^{s}_{\gamma}}^{1/2} \
= \ \Big( \sum_{l, n, \mu} (n + 1)^{s} (n + l + 1)^{s} \, f_{l, n, \mu}^{2} \, h_{l, n}^{2} \Big)^{1/2} \, ,
 \nonumber
\end{align}
where $h_{l, n} \, := \, \| \mcV_{l, \mu}(\bfx) \, P_{n}^{(\gamma, l)}(\rho) \|_{\gamma}$.

\begin{lemma} \label{eqDefDsp}
For $s > 0$, the space $H^{-s}_{\gamma}(\Omega)$ can be characterized as 
\be
\wtilde{H}^{-s}_{\gamma}(\Omega) \ := \ \Big\{ f \, : \, 
f(\bfx) \, = \, \sum_{l, n, \mu} f_{l, n, \mu} \, \mcV_{l, \mu}(\bfx) \, P_{n}^{(\gamma, l)}(\rho) \, , 
\ \mbox{where} \ \sum_{l, n, \mu} (n + 1)^{-s} (n + l + 1)^{-s} \, f_{l, n, \mu}^{2} \, h_{l, n}^{2}  < \infty \Big\} \, . 
\label{altdef1}
\ee
Additionally, for $F \in H^{-s}_{\gamma}(\Omega)$, with representation 
$F  \, = \, \sum_{l, n, \mu} f_{l, n, \mu} \, \mcV_{l, \mu}(\bfx) \, P_{n}^{(\gamma, l)}(\rho) \ := \ f(x)$ 
\be
\| F \|_{H^{-s}_{\gamma}} \ = \ \sup_{g \in H^{s}_{\gamma}(\Omega)} 
\frac{ | F(g) |}{ \| g \|_{H^{s}_{\gamma}(\Omega)}} 
\ = \ \Big( \sum_{l, n, \mu} (n + 1)^{-s} (n + l + 1)^{-s} \, f_{l, n, \mu}^{2} \, h_{l, n}^{2}\Big)^{1/2}  \ := \ \| f \|_{\wtilde{H}^{-s}_{\gamma}} \, .
 \label{fdrmrm}
\ee
\end{lemma}
\textbf{Proof}: To establish the characterization we demonstrate that there is an isometry between 
$H^{-s}_{\gamma}(\Omega)$ and $\wtilde{H}^{-s}_{\gamma}(\Omega)$.

Let $g(\bfx) \, = \, \sum_{l, n, \mu} g_{l, n, \mu} \, \mcV_{l, \mu}(\bfx) \, P_{n}^{(\gamma , l )}(\rho) \in H^{s}_{\gamma}(\Omega)$.

Consider $F(\cdot) \in H^{-s}_{\gamma}(\Omega)$. Then, using the Riesz Representation Theorem,
there exists a unique $k(x) \, = \,\sum_{l, n, \mu} k_{l, n, \mu} \, \mcV_{l, \mu}(\bfx) \, P_{n}^{(\gamma , l )}(\rho)
 \ \in H^{s}_{\gamma}(\Omega)$, 
such that 
\begin{align}
F(g) &= \ \langle k \, ,  \, g \rangle_{H^{s}_{\gamma}} \ = \ 
\sum_{l, n, \mu} (n + 1)^{s} (n + l + 1)^{s} \, k_{l, n, \mu} \, g_{l, n, \mu} \, h_{l, n}^{2}   \nonumber \\
&= \ (f \, ,  \, g )_{\gamma} \, ,   \ \mbox{(see notation in \eqref{deflb})} \nonumber
\end{align}
where $f(\bfx) \, = \, \sum_{l, n, \mu} f_{l, n, \mu} \, \mcV_{l, \mu}(\bfx) \, P_{n}^{(\gamma , l )}(\rho)$,
 for $f_{l, n, \mu} \, = \, (n + 1)^{s} (n + l + 1)^{s} \, k_{l, n, \mu}$.

Note that as
\begin{align*}
\sum_{l, n, \mu}  (n + 1)^{-s} (n + l + 1)^{-s} \, f_{l, n, \mu}^{2} \, h_{l, n}^{2} &= \  
\sum_{l, n, \mu} (n + 1)^{-s} (n + l + 1)^{-s} \, (n + 1)^{2s} (n + l + 1)^{2s} \, k_{l, n, \mu}^{2}  \, h_{l, n}^{2}  \\
 &= \ \sum_{l, n, \mu} (n + 1)^{s} (n + l + 1)^{s} \, k_{l, n, \mu}^{2}  \, h_{l, n}^{2}  \ = \ \| k \|_{H^{s}_{\gamma}}^{2} \, < \infty \, ,
\end{align*}
then $f \in \wtilde{H}^{-s}_{\gamma}(\Omega)$, and from \eqref{fdrmrm}, 
$\| f \|_{\wtilde{H}^{-s}_{\gamma}} \, = \, \| k \|_{H^{s}_{\gamma}}$.

Next, 
\begin{align*}
 \| F \|_{H^{-s}_{\gamma}} &= \ \sup_{g \in H^{s}_{\gamma}(\Omega)}
\frac{ | F(g) |}{ \| g \|_{H^{s}_{\gamma}}}  
\ = \ \sup_{g \in H^{s}_{\gamma}(\Omega)}  \frac{ | ( f \, ,  \, g )_{\gamma} |}{ \| g \|_{H^{s}_{\gamma}}}   \\
&=  \ \sup_{g \in H^{s}_{\gamma}(\Omega)}
 \frac{ | \sum_{l, n, \mu} \, f_{l, n, \mu} \, g_{l, n, \mu} \, h_{l, n}^{2} |}{ \| g \|_{H^{s}_{\gamma}}}   \\
&= \  \sup_{g \in H^{s}_{\gamma}(\Omega)} 
\frac{ | \sum_{l, n, \mu} \, (n + 1)^{-s/2} (n + l + 1)^{-s/2}  \, f_{l, n, \mu} \,(n + 1)^{s/2} (n + l + 1)^{s/2}   g_{l, n, \mu} \,
 h_{l, n}^{2}  |}{ \| g \|_{H^{s}_{\gamma}}}   \\
&\le \ \sup_{g \in H^{s}_{\gamma}(\Omega)}
\frac{ \| f \|_{\wtilde{H}^{-s}_{\gamma}} \, \| g \|_{H^{s}_{\gamma}}}{ \| g \|_{H^{s}_{\gamma}}}  \\
&= \ \| f \|_{\wtilde{H}^{-s}_{\gamma}} \, .
\end{align*}
Additionally, as
\begin{align}
\frac{ | F(k) |}{ \| k \|_{H^{s}_{\gamma}}}  
&= \  \frac{ | ( f \, ,  \, k )_{\gamma} |}{ \| k \|_{H^{s}_{\gamma}}}   
\ = \  \frac{ |  \sum_{l, n, \mu} \, f_{l, n, \mu} \,  k_{l, n, \mu} \, h_{l, n}^{2}  |}{ \| f \|_{\wtilde{H}^{-s}_{\gamma}}}   
\ = \  \frac{ | \sum_{l, n, \mu} \, f_{l, n, \mu} \, \,  (n + 1)^{-s} (n + l + 1)^{-s} \, f_{l, n, \mu} \, h_{l, n}^{2} |}{ \| f \|_{\wtilde{H}^{-s}_{\gamma}}}  
\nonumber  \\ 
&= \ \| f \|_{\wtilde{H}^{-s}_{\gamma}} \, ,  \label{nme3}
\end{align}
then it follows that $ \| F \|_{H^{-s}_{\gamma}} \ = \ \| f \|_{\wtilde{H}^{-s}_{\gamma}}$.

Finally, for $f \in \wtilde{H}^{-s}_{\gamma}(\Omega)$ and any $g \in H^{s}_{\gamma}(\Omega)$,
\begin{align*}
| ( f \, ,  \, g )_{\gamma} |&= \ | \sum_{l, n, \mu} \, f_{l, n, \mu} \,  g_{l, n, \mu} \, h_{l, n}^{2} |  \\
&= \   | \sum_{l, n, \mu} \, (n + 1)^{-s/2} (n + l + 1)^{-s/2} \, f_{l, n, \mu} \, \,  (n + 1)^{s/2} (n + l + 1)^{s/2} \, g_{l, n, \mu} \, h_{l, n}^{2} |  \\
&\le \ \| f \|_{\wtilde{H}^{-s}_{\gamma}} \, \| g \|_{H^{s}_{\gamma}} \, .
\end{align*}
Thus, $f$ defines  a bounded  linear functional on $H^{s}_{\gamma}(\Omega)$, i.e., $f \in H^{-s}_{\gamma}(\Omega)$,
and from \eqref{nme3}, $\| f \|_{H^{-s}_{\gamma}} \, = \,  \| f \|_{\wtilde{H}^{-s}_{\gamma}}$ . \\
\mbox{ } \hfill \qed

\setcounter{equation}{0}
\setcounter{figure}{0}
\setcounter{table}{0}
\setcounter{theorem}{0}
\setcounter{lemma}{0}
\setcounter{corollary}{0}
\setcounter{definition}{0}

\section{Solution of $\big( - \Delta \big)^{\alpot} \, \wtilde{u} \ = \ f$}
\label{sec_ExUng}

In this section we study the fractional Laplace operator as a bounded operator from the
space $\mcH_{\alpot}^{s}(\Omega)$ (defined in equation \eqref{defwtH}) onto $H_{\alpot}^{s}(\Omega)$,
with a bounded inverse.

Recall the definition of the (integral) fractional Laplacian.\\
\textbf{The Fractional Laplacian} \\
In $\real^{d}$, the (integral) fractional Laplacian of a function $u(\bfx)$, 
$\left( - \Delta \right)^{\alpot} u(\bfx)$, is defined as
\be
 \left( - \Delta \right)^{\alpot} u(\bfx) \ := \ \frac{1}{| \gamma_{d}(-\alpha)|} \, \lim_{\eps \rightarrow 0} 
 \int_{\real^{d} \backslash B(\bfx, \eps)} \, \frac{u(\bfx) \, - \, u(\bfy)}{| \bfx \, - \, \bfy |^{d + \alpha}} \, d\bfy \, , 
\label{defFLap}
\ee
where $\gamma_{d}(\alpha) \, := \, 2^{\alpha} \, \pi^{d/2} \, \Gamma(\alpot) / \Gamma( (d - \alpha)/2 )$, and 
$B(\bfx, \eps)$ denotes the ball centered at $\bfx$ with radius~$\eps$.

 Recall $\omega^{\gamma} \, := \, (1 \, - \, r^{2})^{\gamma}$, and let
\be
\mcH^{s}_{\beta}(\Omega) \ := \, \omega^{\beta} \otimes H_{\beta}^{s}(\Omega) \, = \, 
\left\{ \wtilde{v}(\bfx) \, = \, \omega^{\beta} v(\bfx) \, : \, v(\bfx) \in H_{\beta}^{s}(\Omega)  \right\} \, ,  \ \ 
\mbox{with } \ \| \wtilde{v} \|_{\mcH^{s}_{\beta}(\Omega)} \, := \, \| v \|_{H_{\beta}^{s}(\Omega)} \, .
\label{defwtH}
\ee

We have the following property describing the action of the fractional Laplacian on the basis functions for 
$L^{2}_{\alpot}(\Omega)$.
\begin{theorem} \cite[Theorem 3]{dyd171}  \label{genThmLapm}
For $\alpha > 0$, $l, n \, \geq 0$  integer, $\bfx \in \Omega$, 
\begin{align}
   \left( -\Delta \right)^{\frac{\alpha}{2}} (1 \, - \, r^{2})^{\alpot} \, \mcV_{l , \mu}(\bfx) \, P_{n}^{(\alpot ,  l)}(2 r^{2} \, - \, 1) &= \   
   \lambda_{l, n} \, 
  \mcV_{l , \mu}(\bfx) \, P_{n}^{(\alpot ,  l)}(2 r^{2} \, - \, 1) \, ,  \label{bceq0v}  \\
\mbox{where } \ \lambda_{l, n} &= \   2^{\alpha} \, 
  \frac{ \Gamma(n + 1 + \frac{\alpha}{2}) \, \Gamma(n + l + 1 + \frac{\alpha}{2})}%
  {\Gamma(n + 1)\, \Gamma(n + l + 1)} \, .  \nonumber 
\end{align}
\end{theorem}  

Using \eqref{bceq0v} and the function space setting developed in this paper leads to the following theorem.
\begin{theorem} \label{thmon2}
For $s \in \real$, the mapping $\big( - \Delta \big)^{\alpot} \, : \, \mcH_{\alpot}^{s}(\Omega) \longrightarrow 
H_{\alpot}^{s - \alpha}(\Omega)$ is bounded and onto, with a bounded inverse.
\end{theorem}
\textbf{Proof}: Let $\wtilde{u}(\bfx) \ = \ \omega^{\alpot} \, \sum_{l, n, \mu} a_{l, n, \mu} \, \mcV_{l, \mu}(\bfx) \, P_{n}^{(\alpot , l)}(\rho)$.
Then, 
\[
\| \wtilde{u} \|_{\mcH_{\alpot}^{s}}^{2} \ = \ 
  \sum_{l, n, \mu} (n + 1)^{s} \, (n + l + 1)^{s} \, a^{2}_{l, n, \mu} \, \| \mcV_{l, \mu}(\bfx) \, P_{n}^{(\alpot , l)}(\rho) \|^{2}_{L^{2}_{\alpot}}
  \ < \ \infty \, .
\]
Using \eqref{bceq0v},
\[
g(\bfx) \ := \   \big( - \Delta \big)^{\alpot} \, \wtilde{u}(\bfx) 
\ = \ \sum_{l, n, \mu} \lambda_{l, n} \, a_{l, n, \mu} \, \mcV_{l, \mu}(\bfx) \, P_{n}^{(\alpot , l)}(\rho) \, .
\]
Thus,
\begin{align*}
& \| g \|_{H_{\alpot}^{s - \alpha}}^{2} \ = \ 
   \sum_{l, n, \mu} (n + 1)^{s - \alpha} \, (n + l + 1)^{s - \alpha } \,  \lambda_{l, n}^{2} \, a^{2}_{l, n, \mu} \, 
   \| \mcV_{l, \mu}(\bfx) \, P_{n}^{(\alpot , l)}(\rho) \|^{2}_{L^{2}_{\alpot}}  \\
&= \ 
   \sum_{l, n, \mu} (n + 1)^{s - \alpha} \, (n + l + 1)^{s - \alpha } \,
    \Big(  2^{\alpha} \, 
  \frac{ \Gamma(n + 1 + \frac{\alpha}{2}) \, \Gamma(n + l + 1 + \frac{\alpha}{2})}%
  {\Gamma(n + 1)\, \Gamma(n + l + 1)} \Big)^{2} \, 
     a^{2}_{l, n, \mu} \,   \| \mcV_{l, \mu}(\bfx) \, P_{n}^{(\alpot , l)}(\rho) \|^{2}_{L^{2}_{\alpot}}  \\
&\sim \ 
   \sum_{l, n, \mu} (n + 1)^{s - \alpha} \, (n + l + 1)^{s - \alpha } \,
    \big( (n + 1)^{\alpot} \, (n + l + 1)^{\alpot} \big)^{2} \, 
     a^{2}_{l, n, \mu} \,   \| \mcV_{l, \mu}(\bfx) \, P_{n}^{(\alpot , l)}(\rho) \|^{2}_{L^{2}_{\alpot}}  
     \ \ \mbox{(using \eqref{eqStrf})} \\
&= \ 
  \sum_{l, n, \mu} (n + 1)^{s} \, (n + l + 1)^{s} \, a^{2}_{l, n, \mu} \, \| \mcV_{l, \mu}(\bfx) \, P_{n}^{(\alpot , l)}(\rho) \|^{2}_{L^{2}_{\alpot}} \\
&= \   \| \wtilde{u} \|_{\mcH_{\alpot}^{s}}^{2} \, .
\end{align*}
That the mapping is onto follows from the explicit representation of the coefficients of $g(\bfx)$  in terms of the
coefficients of $\wtilde{u}(\bfx)$. \\
\mbox{  } \hfill \qed

The following corollary follows directly from Theorem \ref{thmon2}.
\begin{corollary} \label{corext1}
For $t \in \real$, given $f(\bfx) \in H^{t}_{\alpot}(\Omega)$ there exists a unique $\wtilde{u}(\bfx) \in \mcH^{t + \alpha}_{\alpot}(\Omega)$
satisfying $ \big( - \Delta \big)^{\alpot} \, \wtilde{u}(\bfx) \, = \, f(\bfx) \, , \ \bfx \in \Omega$, 
$\wtilde{u}(\bfx) \, = \, 0 \, , \ \bfx \in \real^{2} \backslash \Omega$, with 
$\| \wtilde{u} \|_{\mcH^{t + \alpha}_{\alpot}(\Omega)} \, \sim \, \| f \|_{H^{t}_{\alpot}(\Omega)}$.
\end{corollary}
\mbox{ } \hfill \qed

\section{Concluding remarks}
In this paper, for $\Omega$ the unit disk in $\real^{2}$, we have shown the equivalence between the weighted 
Sobolev space $H_{\gamma}^{s}(\Omega)$ and the space $\wtilde{H}_{\gamma}^{s}(\Omega)$, defined by the
decay rate of the Fourier-Jacobi coefficients of a function in $L^{2}_{\beta}(\Omega)$. Using
$\wtilde{H}_{\gamma}^{s}(\Omega)$, an explicit representation for $H_{\gamma}^{-s}(\Omega)$ was obtained.
The function space framework developed was then used to give a precise analysis for the fractional Laplace equation
on the unit disk in $\real^{2}$.


\appendix
\setcounter{equation}{0}
\setcounter{figure}{0}
\setcounter{table}{0}
\setcounter{theorem}{0}
\setcounter{lemma}{0}
\setcounter{corollary}{0}
\setcounter{definition}{0}
\section{Proof that $\wtilde{H}^{s}_{\gamma}(\Omega)$, $s \ge 0$, is a family of interpolation spaces}
\label{sec_Intprf}
In this section we present the proof that $\wtilde{H}^{s}_{\gamma}(\Omega)$ is a family of interpolation
spaces with respect to the K-method of interpolation.

\textbf{Proof}:
The proof is similar to an argument used in \cite{sch121}. (See also \cite[Theorem 2.1]{bab011} and \cite[Theorem 4.1]{erv191}.)

Recall for nonnegative integers $k$ and $m$, with $k < s < m$, specifically $s \, = \, (1 - \theta) k \, + \, \theta \, m$ for some $0 < \theta < 1$,
using the $K-$method of interpolation
\[
\| u \|_{X^{s}} \ = \ \left( \int_{0}^{\infty} t^{-2 \theta} \, | K(t, u) |^{2} \, \frac{dt}{t} \right)^{1/2} \, , \ \ 
\mbox{where } \ K(t, u) \ = \ \inf_{v \in X^{m}} \left( \|u - v \|_{X^{k}} \, + \,  t \, \| v \|_{X^{m}} \right) \, .
\]

Noting that for $c$ and $d > 0$, $\frac{1}{\sqrt{2}}(c + d) \, \leq \, (c^{2} + d^{2}) \, \leq \, (c + d)$,  hence we have that
\[
K(t, u) \sim \wtilde{K}(t, u) \ = \ \inf_{v \in X^{m}} \left( \|u - v \|_{X^{k}}^{2} \, + \, t^{2} \, \| v \|_{X^{m}}^{2} \right)^{1/2} \, .
\]

Hence we use
\[
\| u \|_{X^{s}} \ := \ \left( \int_{0}^{\infty} t^{-2 \theta} \, |  \wtilde{K}(t, u) |^{2} \, \frac{dt}{t} \right)^{1/2} \,  .
\]

Let $u(\bfx)$ and $h_{l, n}$ be given by \eqref{defunh}, with $u(\bfx) \in \wtilde{H}^{k}_{\gamma}(\Omega)$ and let 
$v(\bfx) \ = \ \sum_{l, n, \mu} v_{l, n, \mu} \, \mcV(\bfx) \, P_{n}^{(\gamma, l)}(\rho)$ with $v(\bfx) \in \wtilde{H}^{m}_{\gamma}(\Omega)$. Then,
\begin{align*}
 \|u - v \|_{\wtilde{H}^{k}_{\gamma}}^{2} &= \ \sum_{l, n, \mu} (n + 1)^{k} (n + l + 1)^{k} \, (u_{l, n, \mu} - v_{l, n, \mu})^{2} \, 
  h_{l, n}^{2}  \, ,  \\
\mbox{and  }
 \| v \|_{\wtilde{H}^{m}_{\gamma}}^{2} &= \ \sum_{l, n, \mu} (n + 1)^{m} (n + l + 1)^{m} \, v_{l, n, \mu}^{2} \, 
  h_{l, n}^{2}  \, .
\end{align*}

Now,
\begin{align}
\wtilde{K}(t, u)^{2} &= \ \inf_{v \in \wtilde{H}^{m}_{\gamma}(\Omega)} 
 \sum_{l, n, \mu} \left( (n + 1)^{k} (n + l + 1)^{k} \, (u_{l, n, \mu} - v_{l, n, \mu})^{2} \ + \ t^{2} \, 
  (n + 1)^{m} (n + l + 1)^{m} \, v_{l, n, \mu}^{2} \right)  h_{l, n}^{2}  \nonumber \\
&\geq \ \sum_{l, n, \mu}  \inf_{v_{l, n, \mu} \, \in \, \real}  \left( (n + 1)^{k} (n + l + 1)^{k} \, (u_{l, n, \mu} - v_{l, n, \mu})^{2} \ + \ t^{2} \, 
  (n + 1)^{m} (n + l + 1)^{m} \, v_{l, n, \mu}^{2} \right)  h_{l, n}^{2} \, .  \nonumber   
\end{align}

Each of the terms in the summation is minimized for
\be
    v_{l, n, \mu} \ = \ \frac{ u_{l, n, \mu} }{1 \, + \, t^{2} (n + 1)^{m - k} (n + l + 1)^{m - k}} \, .
\label{pol2}
\ee

For this choice of $v_{l, n, \mu}$ we have for $t > 0$,
\begin{align*}
 \| v \|_{\wtilde{H}^{m}_{\gamma}}^{2} 
 &= \ \sum_{l, n, \mu} (n + 1)^{m} (n + l + 1)^{m} \, \frac{ u_{l, n, \mu}^{2} }{(1 \, + \, t^{2} (n + 1)^{m - k} (n + l + 1)^{m - k})^{2}} \, h_{l, n}^{2}   \\
&\leq \ t^{-4}  \sum_{l, n, \mu} \frac{ t^{4} \, (n + 1)^{m} (n + l + 1)^{m} \, u_{l, n, \mu}^{2} }{1 \, + \, t^{2} (n + 1)^{m - k} (n + l + 1)^{m - k}} 
\, h_{l, n}^{2}   \\
&\leq \ t^{-4}  \sum_{l, n, \mu} \, t^{2} \, (n + 1)^{k} (n + l + 1)^{k} \, u_{l, n, \mu}^{2} \, h_{l, n}^{2}   \\
&= \ t^{-2} \, \| u \|_{\wtilde{H}^{k}_{\gamma}}^{2} \, .
\end{align*} 

Hence for the choice of $v_{l, n, \mu}$ given by \eqref{pol2}, $v(\bfx) \in \wtilde{H}^{m}_{\gamma}(\Omega)$. For this choice of  $v_{l, n, \mu}$,
\begin{align*}
& (n + 1)^{k} (n + l + 1)^{k} \, (u_{l, n, \mu} - v_{l, n, \mu})^{2} \ + \ t^{2} \, 
  (n + 1)^{m} (n + l + 1)^{m} \, v_{l, n, \mu}^{2}  \\
& = \   (n + 1)^{k} (n + l + 1)^{k} \, \left( u_{l, n, \mu} \, - \, \frac{ u_{l, n, \mu} }{1 \, + \, t^{2} (n + 1)^{m - k} (n + l + 1)^{m - k}} \right)^{2}  \\
& \quad  \quad  \quad + \ t^{2} \, 
  (n + 1)^{m} (n + l + 1)^{m} \, \left( \frac{ u_{l, n, \mu} }{1 \, + \, t^{2} (n + 1)^{m - k} (n + l + 1)^{m - k}} \right)^{2}  \\
& = \ \frac{t^{2} \,  (n + 1)^{m} (n + l + 1)^{m} }{1 \, + \, t^{2} (n + 1)^{m - k} (n + l + 1)^{m - k}} \, u_{l, n, \mu}^{2} \, .
\end{align*}

Hence,
\[
\wtilde{K}(t, u)^{2} \ = \  \sum_{l, n, \mu} \frac{t^{2} \,  (n + 1)^{m} (n + l + 1)^{m} }%
{1 \, + \, t^{2} (n + 1)^{m - k} (n + l + 1)^{m - k}} \, u_{l, n, \mu}^{2} \, h_{l, n}^{2} \, .
\]

Then,
\begin{align*}
\| u \|^{2}_{\wtilde{H}^{s}_{\gamma}} &= \ 
\int_{0}^{\infty} t^{-2 \theta} \, \wtilde{K}(t, u)^{2} \, \frac{dt}{t}   \\
&= \ \int_{0}^{\infty} t^{-2 \theta} \, \sum_{l, n, \mu} \frac{t^{2} \,  (n + 1)^{m} (n + l + 1)^{m} }%
{1 \, + \, t^{2} (n + 1)^{m - k} (n + l + 1)^{m - k}} \, u_{l, n, \mu}^{2}  \,  h_{l, n}^{2} \, \frac{dt}{t}  \\
&= \  \sum_{l, n, \mu} \, u_{l, n, \mu}^{2}  \,  h_{l, n}^{2} \, \int_{0}^{\infty} t^{-2 \theta} \, \frac{t^{2} \,  (n + 1)^{m} (n + l + 1)^{m} }%
{1 \, + \, t^{2} (n + 1)^{m - k} (n + l + 1)^{m - k}} \, \frac{dt}{t} \, .
\end{align*}

Let $\tau \, = \, t \,  (n + 1)^{(m - k)/2} (n + l + 1)^{(m - k)/2}$. Then $t \, = \, \tau \, (n + 1)^{-(m - k)/2} (n + l + 1)^{-(m - k)/2}$
and $dt \, = \,  (n + 1)^{-(m - k)/2} (n + l + 1)^{-(m - k)/2} \, d\tau$. Substituting,
\begin{align*}
\| u \|^{2}_{\wtilde{H}_{\gamma}^{s}} &= \ \sum_{l, n, \mu} \, u_{l, n, \mu}^{2}  \,  h_{l, n}^{2} \, \int_{0}^{\infty} 
\tau^{-2 \theta} \, (n + 1)^{(m - k) \theta} (n + l + 1)^{(m - k) \theta} \, 
\frac{\tau^{2} \,  (n + 1)^{k} (n + l + 1)^{k} }%
{1 \, + \, \tau^{2}} \, \frac{d\tau}{\tau} \\
&= \ \sum_{l, n, \mu} \,  (n + 1)^{k \, + \, (m - k) \theta} (n + l + 1)^{k \, + \, (m - k) \theta} \,  u_{l, n, \mu}^{2}  \,  h_{l, n}^{2} 
 \, \int_{0}^{\infty} \frac{\tau^{1 - 2\theta}}{1 + \tau^{2}} \, d\tau \\
&= \ C_{\theta}^{2} \, \sum_{l, n, \mu} \,  (n + 1)^{k \, + \, (m - k) \theta} (n + l + 1)^{k \, + \, (m - k) \theta} \,  u_{l, n, \mu}^{2}  \,  h_{l, n}^{2}  \\
&= \ C_{\theta}^{2} \, \sum_{l, n, \mu} \,  (n + 1)^{s} (n + l + 1)^{s} \,  u_{l, n, \mu}^{2}  \,  h_{l, n}^{2}  \, ,
\end{align*}
where $C_{\theta}^{2} \, = \, \int_{0}^{\infty} \frac{\tau^{1 - 2\theta}}{1 + \tau^{2}} \, d\tau  \, < \, \infty$. 

Hence the family of spaces $\wtilde{H}_{\gamma}^{s}(\Omega)$ represent 
interpolation spaces with respect to the K-method of interpolation. \\
\mbox{ } \hfill \qed

\setcounter{equation}{0}
\setcounter{figure}{0}
\setcounter{table}{0}
\setcounter{theorem}{0}
\setcounter{lemma}{0}
\setcounter{corollary}{0}
\setcounter{definition}{0}
\section{Explicit representation of $| u |_{H^{k}_{\gamma}(\Omega)} $ }
\label{sec_snorm2}
In this section we obtain an explicit expression for $| u |_{H^{k}_{\gamma}(\Omega)} $, where
\be
 u(\bfx) \ = \ \sum_{l, n, \mu} a_{l, n, \mu} \, \mcV_{l, \mu}(\bfx) \, P_{n}^{(\gamma , l )}(\rho) \, ,
 \ \ ( \mbox{recall } \rho \, = \, 2r^{2} - 1) \, , 
\label{repApu}
\ee
in terms of the coefficients $a_{l, n, \mu}$ and $\| \mcV_{l, \mu}(\bfx) \, P_{n}^{(\gamma , l )}(\rho) \|_{\gamma}$. 

In Section \ref{ssec_null} we determine those basis functions in 
$\left\{ \mcV_{l, \mu}(\bfx) P_{n}^{(\gamma , l )}(\rho) \right\}_{l, n, \mu}$ whose all $k^{th}$ order derivatives
are zero.

In Section \ref{ssec_nexp} several equations for the basis functions in different norms are calculated, which
are useful in determining expressions for the semi norms.

\textbf{Remark}: For convenience in the analysis below, if $l < 0$ or $n < 0$, then $a_{l, n, \mu} = 0$. In addition, if
$l = 0$ and $\mu = -1$, then $a_{l, n, \mu} = 0$.

\subsection{Basis functions whose all $k^{th}$ order derivatives are zero}
\label{ssec_null}

In this section we determine those basis functions in 
$\left\{ \mcV_{l, \mu}(\bfx) P_{n}^{(\gamma , l )}(\rho) \right\}_{l, n, \mu}$ whose all $k^{th}$ order derivatives
are zero.

\begin{lemma} \label{lmaker1}
We have the following: \\
(i) $\mcV_{l, \mu}(\bfx)$ is comprised of terms $x^{a} \, y^{b}$, $a, \, b \in \Natural_{0}$, of degree exactly $l$, \\
(ii) $\mcV_{l, 1}(\bfx)$ contains the term $x^{l}$, \\
(iii) $\mcV_{l, -1}(\bfx)$ contains the term $l \, y \, x^{l-1}$, and does not contain the term $x^{l}$.
\end{lemma}
\textbf{Proof}: The proof is by induction. Note that
\[
\begin{array}{rlcrl}
\mcV_{1, 1}(\bfx) &= \ r \cos( \varphi ) \ = \ x \, , &\quad  &\mcV_{1, 1}(x) &= \ r \sin( \varphi ) \ = \ y \\
\mcV_{2, 1}(\bfx) &= \ r^{2} \cos( 2 \varphi ) &\quad  &\mcV_{2, -1}(x) &= \ r^{2} \sin( 2 \varphi )  \\
                        &= r^{2} ( \cos(\varphi)^{2} \, - \, \sin(\varphi)^{2}  &\quad   &    &= r^{2} \, 2 \sin( \varphi ) \cos( \varphi )  \\
                        &= x^{2} \, - \, y^{2}      & \quad  &  &= 2 \, x  \, y \, .
\end{array}                        
\]
\underline{Induction step for $\mcV_{l, 1}(\bfx)$} 
\begin{align*}
\mcV_{l, 1}(\bfx) &= \ r^{l} \cos(l \varphi)   \ = \ r^{l} \, \cos( \varphi \, + \, (l - 1) \varphi )  \\
&= \ r^{l} \left( \cos(\varphi) \, \cos((l - 1) \varphi) \ - \  \sin(\varphi) \, \sin((l - 1) \varphi)  \right)  \\
&= \ x \, \mcV_{l-1 \, , 1}(x) \ - \ y \, \mcV_{l , -1}(x) \, .
\end{align*}
By the induction assumption, $\mcV_{l-1 \, , 1}(\bfx)$ and $\mcV_{l-1 \, , -1}(\bfx)$ only contain polynomial terms of degree $(l - 1)$.
Hence $\mcV_{l, 1}(\bfx)$ only contains terms of degree $l$. Moreover, as $V_{l-1 \, , 1}(\bfx)$ contains the term $x^{l - 1}$ then
$\mcV_{l, 1}(\bfx)$ contains the term $x^{l}$.

\underline{Induction step for $\mcV_{l, -1}(\bfx)$} 
\begin{align*}
\mcV_{l, -1}(\bfx) &= \ r^{l} \sin(l \varphi)   \ = \ r^{l} \, \sin( \varphi \, + \, (l - 1) \varphi )  \\
&= \ r^{l} \left( \sin(\varphi) \, \cos((l - 1) \varphi) \ + \  \cos(\varphi) \, \sin((l - 1) \varphi)  \right)  \\
&= \ y \, \mcV_{l-1 \, , 1}(x) \ + \ x \, \mcV_{l , -1}(x) \, .
\end{align*}
Similar to $\mcV_{l, 1}(\bfx)$, it follows that $\mcV_{l, -1}(\bfx)$ only contains terms of degree $l$. \\
Using the induction assumption that $\mcV_{l-1 \, , 1}(\bfx)$ contains $x^{l - 1}$ and $\mcV_{l-1 \, , -1}(\bfx)$ contains $(l - 1) \, x^{l - 2} \, y$,
then $\mcV_{l , -1}(\bfx)$ contains the term $l \, x^{l - 1} \, y$. As $\mcV_{l-1 \, , -1}(\bfx)$ does not contains $x^{l - 1}$, then 
$\mcV_{l  , 1}(\bfx)$ does not contains $x^{l}$. \\
\mbox{  } \hfill \qed

\begin{lemma}  \label{lmaker2}
$\mcV_{l , \mu}(\bfx) \, P_{n}^{(\gamma , l)}(\rho)$ is a nonzero polynomial in $x$ and $y$ of degree exactly $2n \, + \, l$.
\end{lemma}
\textbf{Proof}: Recall that $P_{n}^{(\gamma , l)}(t)$ is a polynomial of degree $n$ in $t$. Hence $P_{n}^{(\gamma , l)}(2r^{2} - 1)$
is a polynomial of degree $2n$ in $x$ and $y$. Additionally $P_{n}^{(\gamma , l)}(2r^{2} - 1)$ only contains even powers of 
$x$ and $y$. In particular, $P_{n}^{(\gamma , l)}(2r^{2} - 1)$ contains the term $c \, x^{2n}$ for $c \neq 0$.
Hence, $\mcV_{l , 1}(\bfx) \, P_{n}^{(\gamma , l)}(2r^{2} - 1)$ is a polynomial of degree exactly $(2n \, + \, l)$, as it contains the
term $c \, x^{2n + l}$. Also, $\mcV_{l , -1}(\bfx) \, P_{n}^{(\gamma , l)}(2r^{2} - 1)$ is a polynomial of degree exactly $(2n \, + \, l)$, 
as it contains the term $c \, l \, y \, x^{2n + l - 1}$. For $\mcV_{0, 1}(\bfx) \, = \, \frac{1}{2}$ the result follows for the second statement 
in the proof. \\
\mbox{  } \hfill \qed

\begin{lemma}  \label{lmaker3}
Let $\mcG_{k} \, : \ X \, := \, \left\{ \mcV_{l , \mu}(\bfx) \, P_{n}^{(\gamma , l)}(2r^{2} - 1) \right\}_{l, n, \mu} \, \longrightarrow \, 
\left( C(\Omega) \right)^{k + 1}$, be defined by
\[
  \mcG_{k} f(\bfx) \ := \ 
  \left[ \frac{\partial^{k} f(\bfx)}{\partial x^{k}} , \ldots ,  \frac{\partial^{k} f(\bfx)}{\partial y^{j} \partial x^{k-j}}, \ldots ,
  \frac{\partial^{k} f(\bfx)}{\partial y^{k}} \right] \, .
\]
Additionally, let
\[
 Null(  \mcG_{k} ) \ := \ \left\{ f(\bfx) \in X \, : \, \mcG_{k} f(\bfx) \ = \ [0, 0, \ldots, 0] \right\} \, .
\]
Then,
\[
 Null(  \mcG_{k} ) \ = \ \cup_{2n + l < k} \left\{ \mcV_{l , \mu}(\bfx) \, P_{n}^{(\gamma , l)}(2r^{2} - 1) \right\} \, .
\]
\end{lemma}
\textbf{Proof}: Note that
$\mcV_{l , \mu}(\bfx) \, P_{n}^{(\gamma , l)}(2r^{2} - 1) \in  Null(  \mcG_{k} )$ if and only if all its $k^{th}$ order derivatives are 0.
From Lemma \ref{lmaker2}, $\mcV_{l , \mu}(\bfx) \, P_{n}^{(\gamma , l)}(2r^{2} - 1)$ is a polynomial in $x$ and $y$ of degree exactly 
$2n \, + \, l$. Hence, all the $k^{th}$ order derivatives of $\mcV_{l , \mu}(\bfx) \, P_{n}^{(\gamma , l)}(2r^{2} - 1)$ are 0 if and 
only if $2n \, + \, l < k$.  \\
\mbox{ } \hfill \qed

\subsection{Some useful norm expressions for the basis functions}
\label{ssec_nexp}

In this section we present a number of equations for the basis functions
measured in different norms.

Let $C_{l, \mu} \ = \ \left\{ \begin{array}{l}
                                  \pi/2 , \ \mbox{if } l = 0, \ \mu = 1,  \\
                                  \pi , \ \mbox{for } l \ge 1, 
                                  \end{array}  \right. $
\begin{lemma} \label{lmaeq3}
For $l, \, n, \, k \ge 0$, $(l, \, \mu) \neq (0, -1)$,
\be 
\frac{\| \mcV_{l, \mu}(\bfx) \, P_{n}^{(\gamma + k \, , \, l)}(\rho) \|^{2}_{L^{2}_{\gamma + k}}}%
{  \| \mcV_{l, \mu}(\bfx) \, P_{n}^{(\gamma , l)}(\rho) \|^{2}_{L^{2}_{\gamma}} }
  \ = \ \frac{2n + \gamma + l + 1}{2n + \gamma + l + 1 + k} \, 
 \frac{ \Pi_{s = 1}^{k} \, (n + \gamma + s)}{ \Pi_{s = 1}^{k} \, (n + \gamma + l + s)} \,  .
\label{poit1}
\ee
\end{lemma}
\textbf{Proof}:
\begin{align*}
& \|  \mcV_{l  , \mu}(\bfx) \, P_{n}^{(\gamma + k \,  , \,   l)}(\rho)  \|^{2}_{L^{2}_{\gamma + k}} 
\ = \ \int_{\Omega} \omega^{\gamma + k} \left( \mcV_{l  , \mu}(\bfx) \, P_{n}^{(\gamma + k \, ,  \, l)}(2 r^{2} - 1)  \right)^{2} \, d\Omega  \\
&= \  C_{l, \mu}  \,  \int_{t = -1}^{1} \left( \frac{1 - t}{2} \right)^{\gamma + k} \, \left( \frac{1 + t}{2} \right)^{l} \,   
   \left( P_{n}^{(\gamma + k \,  ,  \,  l)}(t)  \right)^{2} \, \frac{1}{4} \, dt \, , \ \ \ \mbox{ for } t \, = \, 2 r^{2} - 1 \, ,   \\
&= \  C_{l, \mu}  \, 2^{-(\gamma + k + l + 2)} \, |\| P_{n}^{(\gamma + k \, ,  \, l)}(t)  |\|^{2}    \\
&= \  C_{l, \mu}  \, 2^{-(\gamma + k + l + 2)} \, \frac{2^{(\gamma + k + l + 1)}}{2n + \gamma + k + l + 1} \, 
\frac{\Gamma(n + \gamma + k + 1) \, \Gamma(n + l + 1)}{\Gamma(n + 1) \, \Gamma(n + \gamma + k + l + 1)} \,
 \ \mbox{ using } \eqref{spm22}  \\
&= \ C_{l, \mu}  \, 2^{-(\gamma + l + 2)} \, \frac{2n + \gamma + l + 1}{2n + \gamma + k + l + 1} \, 
\frac{\Gamma(n + \gamma + k + 1)}{ \Gamma(n + \gamma + 1)} \,
\frac{\Gamma(n + \gamma + l + 1)}{\Gamma(n + \gamma + k + l + 1)} \, |\| P_{n}^{(\gamma ,  l)}(t)  |\|^{2} \\
&= \ \frac{2n + \gamma + l + 1}{2n + \gamma + k + l + 1} 
 \frac{ \Pi_{s = 1}^{k} \, (n + \gamma + s)}{ \Pi_{s = 1}^{k} \, (n + \gamma + l + s)} \, 
 \| \mcV_{l, \mu}(\bfx) \, P_{n}^{(\gamma , l)}(\rho) \|^{2}_{L^{2}_{\gamma}} \, .
 \end{align*}
\mbox{ } \hfill \qed

\begin{lemma} \label{lmaeq4}  
For $l, \, n, \, j, \, m \geq 0$, $(l, \, \mu) \neq (0, -1)$,
\begin{align}
& \frac{ \| \mcV_{l+j \, , \, \mu}(\bfx) \, P_{n+m}^{(\gamma , \, l + j )}(\rho) \|^{2}_{L^{2}_{\gamma}} }%
{ \| \mcV_{l, \mu}(\bfx) \, P_{n}^{(\gamma , l)}(\rho) \|^{2}_{L^{2}_{\gamma}} }  \nonumber \\
& \quad \quad \quad \ = \ \frac{C_{l+j , \mu}}{C_{l , \mu}} \, \frac{2n + \gamma + l + 1}{2n + 2m + \gamma + l + j + 1} \, 
 \frac{ \Pi_{s = 1}^{m} \, (n + \gamma + s)}{ \Pi_{s = 1}^{m} \, (n +  s)} \, 
 \frac{ \Pi_{s = 1}^{m+j} \, (n + l + s)}{ \Pi_{s = 1}^{m+j} \, (n + \gamma + l + s)} \, 
 \, .
\label{poit11}
\end{align}
\end{lemma}
\textbf{Proof}:
\begin{align*}
& \|  \mcV_{l+j \, , \, \mu}(\bfx) \, P_{n+m}^{(\gamma   , \,   l + j )}(\rho)  \|^{2}_{L^{2}_{\gamma}} 
\ = \ \int_{\Omega} \omega^{\gamma } \left( \mcV_{l+j \, , \, \mu}(\bfx) \, P_{n+m}^{(\gamma   , \,   l + j )}(2 r^{2} - 1)  \right)^{2} \, d\Omega  \\
&= \  C_{l+j, \mu}  \,  \int_{t = -1}^{1} \left( \frac{1 - t}{2} \right)^{\gamma} \, \left( \frac{1 + t}{2} \right)^{l + j} \,   
   \left( P_{n+m}^{(\gamma  ,  \,  l + j)}(t)  \right)^{2} \, \frac{1}{4} \, dt \, , \ \ \ \mbox{ for } t \, = \, 2 r^{2} - 1 \, ,   \\
&= \  C_{l+j, \mu}  \, 2^{-(\gamma + l + j + 2)} \, |\| P_{n+m}^{(\gamma ,  \, l + j)}(t)  |\|^{2}    \\
&= \  C_{l+j, \mu}  \, 2^{-(\gamma + l + j + 2)} \, \frac{2^{(\gamma + l + j + 1)}}{2n + 2m + \gamma + l + j + 1} \, 
\frac{\Gamma(n + m + \gamma + 1) \, \Gamma(n + m + l + j + 1)}{\Gamma(n + m + 1) \, \Gamma(n + m + \gamma + l + j + 1)} \,
 \ \mbox{ using } \eqref{spm22}  \\
&= \ C_{l+j, \mu}  \, 2^{-(\gamma + l + 2)} \, \frac{2n + \gamma + l + 1}{2n + 2m + \gamma + l + j + 1} \, 
\frac{\Gamma(n + m + \gamma + 1)}{\Gamma(n + \gamma + 1)} \,
\frac{\Gamma(n + 1)}{\Gamma(n + m + 1)} \,   \\
& \quad \quad  \quad \quad \cdot \ 
\frac{\Gamma(n +  m + l + j + 1)}{\Gamma(n + l + 1)} \,
\frac{\Gamma(n +  \gamma + l + 1)}{\Gamma(n + m + \gamma + l + j + 1)} \,
 |\| P_{n}^{(\gamma ,  l)}(t)  |\|^{2} \\
&= \ \frac{C_{l+j, \mu}}{C_{l, \mu}}   \,    \frac{2n + \gamma + l + 1}{2n + 2m + \gamma + l + j + 1}  
 \frac{ \Pi_{s = 1}^{m} \, (n + \gamma + s)}{ \Pi_{s = 1}^{m} \, (n + s)} \, 
 \frac{ \Pi_{s = 1}^{m + j} \, (n + l + s)}{ \Pi_{s = 1}^{m + j} \, (n + \gamma + l + s)} \, 
 \| \mcV_{l, \mu}(\bfx) \, P_{n}^{(\gamma , l)}(\rho) \|^{2}_{L^{2}_{\gamma}} \, .
 \end{align*}
\mbox{ } \hfill \qed

\begin{lemma} \label{lmaeq5}
For $l, \, n \ge 0$ and $m \geq j \geq 0$ with $(l - j , \, \mu) \neq (0, -1)$,
\begin{align} 
& \frac{ \| \mcV_{l-j \, , \, \mu}(\bfx) \, P_{n+m}^{(\gamma , \, l - j )}(\rho) \|^{2}_{L^{2}_{\gamma}} }%
{ \| \mcV_{l, \mu}(\bfx) \, P_{n}^{(\gamma , l)}(\rho) \|^{2}_{L^{2}_{\gamma}} }  \nonumber \\
& \quad \quad \quad \ = \ \frac{C_{l-j , \mu}}{C_{l , \mu}} \, \frac{2n + \gamma + l + 1}{2n + 2m + \gamma + l - j + 1} \, 
 \frac{ \Pi_{s = 1}^{m} \, (n + \gamma + s)}{ \Pi_{s = 1}^{m} \, (n +  s)} \, 
 \frac{ \Pi_{s = 1}^{m-j} \, (n + l + s)}{ \Pi_{s = 1}^{m-j} \, (n + \gamma + l + s)} \, 
 \, .
\label{poit21}
\end{align}
\end{lemma}
\textbf{Proof}: The proof follows completely analogous to that of Lemma \ref{lmaeq4}.\\
\mbox{ } \hfill \qed

\subsection{Explicit representation of  $| u |_{H^{2}_{\gamma}} $ }
\label{ssec_H2}
In this section we obtain an explicit representation for $| u |_{H^{2}_{\gamma}} $.

\begin{lemma}  \label{thmsnm2}
We have that for $u(\bfx)$ given by \eqref{repApu},
\be
| u |_{H^{2}_{\gamma}}^{2} \ \sim \  \sum_{2n + l \ge 2, \mu} (n + 1)^{2} \, (n + l + 1)^{2} \, a_{l, n, \mu}^{2} \, 
  \| \mcV_{l, \mu}(x) P_{n}^{(\gamma , l)}(\rho) \|_{L^{2}_{\gamma}}^{2} \, . 
\label{repsn2}
\ee
\end{lemma}
\textbf{Proof}: 
Note that 
\begin{align*}
| u |_{H^{2}_{\gamma}}^{2} &= \ 
\left\| \frac{\partial}{\partial x} \left( \frac{\partial u}{\partial x} \right) \right\|^{2}_{L^{2}_{\gamma + 2}} \, + \, 
\left\| \frac{\partial}{\partial y} \left( \frac{\partial u}{\partial x} \right) \right\|^{2}_{L^{2}_{\gamma + 2}} \, + \, 
\left\| \frac{\partial}{\partial x} \left( \frac{\partial u}{\partial y} \right) \right\|^{2}_{L^{2}_{\gamma + 2}} \, + \, 
\left\| \frac{\partial}{\partial y} \left( \frac{\partial u}{\partial y} \right) \right\|^{2}_{L^{2}_{\gamma + 2}}   \\
&= \ 
\left|  \frac{\partial u}{\partial x}  \right|_{H^{1}_{\gamma + 1}}^{2} \, + \, 
\left|  \frac{\partial u}{\partial y}  \right|_{H^{1}_{\gamma + 1}}^{2} \, . 
\end{align*}
So, we proceed in a similar fashion to the proof of Lemma \ref{snH01} for the functions 
$ \frac{\partial u}{\partial x}$ and $ \frac{\partial u}{\partial y}$.

From  \eqref{poit2p5} and \eqref{poit2}, (with $\gamma \rightarrow \, \gamma + 1$)
\begin{align}
\frac{\partial}{\partial x} \left( \frac{\partial u}{\partial x} \right)  &= \ 
\sum_{l, n, \mu} \, _{x x}A_{l, n, \mu} \, \mcV_{l , \mu}(x) P_{n}^{(\gamma + 2 \, , \, l)}(\rho) \, ,  \label{ghg10} \\
\mbox{ where } \ _{x x}A_{l, n, \mu} &= \  _{x}A_{l+1 , n, \mu} (n + l + 1) \, + \, _{x}A_{l-1 , n+1 , \mu} (n + \gamma + l + 2) \, . \nonumber 
\end{align}
Also, using \eqref{poit2p5} and \eqref{poit4}, (with $\gamma \rightarrow \, \gamma + 1$)
\begin{align}
\frac{\partial}{\partial y} \left( \frac{\partial u}{\partial x} \right)  &= \ 
\sum_{l, n, \mu} \, _{y x}A_{l, n, \mu} \, \mcV_{l , \mu^{*}}(x) P_{n}^{(\gamma + 2 \, , \, l)}(\rho) \, ,   \label{ghg11} \\
\mbox{ where } \ _{y x}A_{l, n, \mu} &= \   (\mp) _{x}A_{l+1 , n, \mu} (n + l + 1) \, + \, (\pm) _{x}A_{l-1 , n+1 , \mu} (n + \gamma + l + 2) \, .
  \label{ghg12}
\end{align}

Analogous to \eqref{poit3}, \eqref{poit5p5} and \eqref{poit8p5} we obtain,
\begin{align}
& \left|  \frac{\partial u}{\partial x}  \right|_{H^{1}_{\gamma + 1}}^{2}
\ = \ \sum_{n} \left(  _{x x}A_{0, n, 1}^{2}  \, + \, _{y x}A_{0, n, -1}^{2}  \right)
 \, \| \mcV_{0 , 1}(\bfx) P_{n}^{(\gamma + 2 \, , \, 0)}(\rho) \|_{\gamma + 2}^{2} \, ,   \nonumber   \\
& \quad \hspace{1.0in} +  \sum_{l \ge 1, \, n, \mu} \left(  _{x x}A_{l, n, \mu}^{2}  \, + \, _{y x}A_{l, n, \mu}^{2}  \right)
 \, \| \mcV_{l , 1}(\bfx) P_{n}^{(\gamma + 2 \, , \, l)}(\rho) \|_{\gamma + 2}^{2} \, ,  \nonumber  \\
%
%
%
%
%
& \ = \  \sum_{n}  \big( \, _{x}A_{1 , n, 1}^{2} \ + \ _{x}A_{1 , n, -1}^{2} \big) \, (n + 1)^{2}  
 \, \| \mcV_{0 , 1}(\bfx) P_{n}^{(\gamma + 2 \, , \, 0)}(\rho) \|^{2}_{\gamma + 2}   \nonumber \\
 & \ \ \ + \ 2 \,  \sum_{l \ge 1, \, n, \mu}  
 \Big(  \, _{x}A_{l+1 , n, \mu}^{2} (n + l + 1)^{2} \, + \, _{x}A_{l-1 , n+1, \mu}^{2} (n + \gamma + l + 2)^{2} \Big)
 \, \| \mcV_{l , 1}(\bfx) P_{n}^{(\gamma + 2 \, , \, l)}(\rho) \|^{2}_{\gamma + 2}   \label{mjh0p5} \\
%
%
%
%
& \ \sim \ \sum_{2n + l \ge 1 \, , \mu} (n + 1) (n + l + 1) \,  _{x}A_{l , n, \mu}^{2} \, 
\| \mcV_{l, \mu}(x) P_{n}^{(\gamma + 1 \,  , l)}(\rho) \|^{2}_{\gamma + 1}  \, .
\label{mjh1}
\end{align}

Similarly, 
\begin{align}
& \left|  \frac{\partial u}{\partial y}  \right|_{H^{1}_{\gamma + 1}}^{2}
\ = \ \sum_{n} \left( _{x y}A_{0, n, 1}^{2}  \, + \, _{y y}A_{0, n, -1}^{2}  \right)
 \, \| \mcV_{0 , 1}(\bfx) P_{n}^{(\gamma + 2 \, , \, 0)}(\rho) \|_{\gamma + 2}^{2} \, ,   \nonumber   \\
&  \quad \hspace{1.0in}  +  \sum_{l \ge 1, \, n, \mu} \left( _{x y}A_{l, n, \mu}^{2}  \, + \, _{y y}A_{l, n, \mu}^{2}  \right)
 \, \| \mcV_{l , 1}(\bfx) P_{n}^{(\gamma + 2 \, , \, l)}(\rho) \|_{\gamma + 2}^{2} \, ,  \nonumber  \\
& \ = \  \sum_{n}  \big( \, _{y}A_{1 , n, 1}^{2} \ + \ _{y}A_{1 , n, -1}^{2} \big) \, (n + 1)^{2}  
 \, \| \mcV_{0 , 1}(\bfx) P_{n}^{(\gamma + 2 \, , \, 0)}(\rho) \|^{2}_{\gamma + 2}   \nonumber \\
 & \ \ + \ 2 \,  \sum_{l \ge 1, \, n, \mu}  
 \Big( \, _{y}A_{l+1 , n, \mu}^{2} (n + l + 1)^{2} \, + \, _{y}A_{l-1 , n+1, \mu}^{2} (n + \gamma + l + 2)^{2} \Big)
 \, \| \mcV_{l , 1}(\bfx) P_{n}^{(\gamma + 2 \, , \, l)}(\rho) \|^{2}_{\gamma + 2}   \label{mjh1p5} \\
& \ \sim \ \sum_{2n + l \ge 1 \, , \mu} (n + 1) (n + l + 1) \,  _{y}A_{l , n, \mu}^{2} \, 
\| \mcV_{l, \mu}(x) P_{n}^{(\gamma + 1 \,  , l)}(\rho) \|^{2}_{\gamma + 1}  \, .
\label{mjh2}
\end{align}


Combining \eqref{mjh1} and \eqref{mjh2}, 
\begin{align}
| u |_{H^{2}_{\gamma}}^{2} &= \ 
 \left|  \frac{\partial u}{\partial x}  \right|_{H^{1}_{\gamma + 1}}^{2} \, + \, 
 \left|  \frac{\partial u}{\partial y}  \right|_{H^{1}_{\gamma + 1}}^{2}     \nonumber  \\
&= \ \sum_{l, n, \mu} \Big(  \, _{xx}A_{l , n, \mu}^{2} \, + \, _{yx}A_{l , n, \mu}^{2} \, + \,
 _{xy}A_{l , n, \mu}^{2} \, + \, _{yy}A_{l , n, \mu}^{2} \Big) \, 
  \| \mcV_{l, \mu}(\bfx) P_{n}^{(\gamma + 2 \,  , l)}(\rho) \|^{2}_{\gamma + 2}   \label{mjh25}  \\
&\sim \ \sum_{2n + l \ge 1 \, , \mu} (n + 1) (n + l + 1) \, ( _{x}A_{l , n, \mu}^{2} \, + \, _{y}A_{l , n, \mu}^{2} )
\| \mcV_{l, \mu}(\bfx) P_{n}^{(\gamma + 1 \,  , l)}(\rho) \|^{2}_{\gamma + 1}  \,  . \nonumber 
\end{align}
Following the steps from \eqref{poit8p3} to \eqref{poit8p5} 
(i.e., (i) replace $( _{x}A_{l , n, \mu}^{2} \, + \, _{y}A_{l , n, \mu}^{2} )$ by 
$2 \, \big( (n \, + \, l \, + \, 1)^{2} \, a_{l+1, n, \mu}^{2} \ + \ (n \, + \, \gamma \, + \, l \, + \, 1)^{2} \, a_{l-1, n+1, \mu}^{2} \big)$, 
(ii) expand the resulting product with $\| \mcV_{l, \mu}(\bfx) P_{n}^{(\gamma + 1 \,  , l)}(\rho) \|^{2}_{\gamma + 1} $,
(iii) use Lemmas \ref{lmaeq3}, \ref{lmaeq4} and \ref{lmaeq5} to replace 
$\| \mcV_{l, \mu}(\bfx) P_{n}^{(\gamma + 1 \,  , l)}(\rho) \|^{2}_{\gamma + 1} $ with corresponding terms which match the
subscripts of the $a_{., ., .}$ terms, (iv) reindex the terms so that the subscript on the $a_{., ., .}$ terms is $l, n, \mu$ ) 
yields
\begin{align}
| u |_{H^{2}_{\gamma}}^{2}
%
&\sim \ \sum_{2n + l \ge 2 \, , \mu} (n + 1)^2 (n + l + 1)^2 \, a_{l, n, \mu}^{2} \, 
\| \mcV_{l, \mu}(\bfx) P_{n}^{(\gamma  , l)}(\rho) \|^{2}_{L^{2}_{\gamma}}  \, ,  \label{mjh4} 
\end{align}
where in the last step we have used Lemma \ref{lmaker3}.  \\
\mbox{ } \hfill \qed

Some additional notation is need for the proof of Lemma \ref{snH1}. \\
Let $\bfX^{k} \, := \, \{x , y \} \times \{x , y \} \ldots  \times \{x , y \} \, = \, \{x , y \}^{k}$, and for $\bfp \in \bfX^{k}$, say
$\bfp \, = \, [x, x, y, x, \ldots, x]$ we define
\[
  \partial_{\bfp} u \, := \, \frac{\partial^{k} u}{\partial x \, \partial x \, \partial y \, \partial x \, \ldots \partial x} \, .
\]
Note that the set of all $k^{th}$ order partial derivatives of $u$ can be represented as $\{ \partial_{\bfp} u \}_{\bfp \in \bfX^{k}}$,
and the set of $(k + 1)^{st}$ order derivatives can be expressed as
\[
   \left\{ \frac{\partial}{\partial x} \left( \partial_{\bfp} u \right) \right\}_{\bfp \in \bfX^{k}} \cup
    \left\{ \frac{\partial}{\partial y} \left( \partial_{\bfp} u \right) \right\}_{\bfp \in \bfX^{k}}  \, .
\]    

Similarly, for $\bfp = [x] \in \bfX^{1}$, we use $_{\bfp}A_{l, n, \mu} \, = \, _{x}A_{l, n, \mu}$ (see \eqref{poit2p5}),
for $\bfp = [y, x] \in \bfX^{2}$, we use $_{\bfp}A_{l, n, \mu} \, = \, _{y x}A_{l, n, \mu}$ (see \eqref{ghg12}), etc.  \\

\textbf{Proof of Lemma \ref{snH1}} \\
The proof is by induction. We assume the result is true for $s = k$ and show that it then holds for $s \, = \, k + 1$, i.e, we assume:
for $\bfp \in \bfX^{k}$
\be
  \partial_{\bfp} u(\bfx) \ = \ \sum_{l, n, \mu} \, _{\bfp}A_{l, n, \mu} \, \mcV_{l, \sigma}(\bfx) \, 
  P_{n}^{(\gamma + k \, , \, l)}(\rho) \ \ \mbox{(see \eqref{poit2p5}, \eqref{poit5}, \eqref{ghg10}, \eqref{ghg11})} \, , 
  \label{mkl1}
 \ee
where $\sigma$ is either $\mu$ or $\mu^{*}$, and
\begin{align}
| u |^{2}_{H^{k}_{\gamma}} 
&= \ \sum_{l, n, \mu}  \sum_{\bfp \in \bfX^{k}} \, _{\bfp}A_{l, n, \mu}^{2} \, \| \mcV_{l, \mu}(\bfx) \,
 P_{n}^{(\gamma + k \, , \, l)}(\rho) \|^{2}_{L^{2}_{\gamma + k}} \ \ \mbox{(see \eqref{poit8p3}, \eqref{mjh25})}  \label{mkl2} \\
&\sim \ \sum_{2n + l \ge k \, , \mu} (n + 1)^k (n + l + 1)^k \, a_{l, n, \mu}^{2} \, 
\| \mcV_{l, \mu}(x) P_{n}^{(\gamma  , l)}(\rho) \|^{2}_{L^{2}_{\gamma}} \ \ \mbox{(see \eqref{poit8p5}, \eqref{mjh4})} 
 \, .   \label{mkl3}
\end{align}

Proceeding as in the cases for $k = 1$ (Lemma \ref{snH1}) and $k = 2$ (Lemma \ref{thmsnm2}), for each of the $k^{th}$ order
partial derivative terms \eqref{mkl1} we compute their partial derivative with respect to $x$ and $y$ to obtain
(using Corollary \ref{gradwR2})
\begin{align*}
\frac{\partial}{\partial x} \partial_{\bfp} u(\bfx) &= \ 
\sum_{l, n, \mu} \, _{[x , \bfp]}A_{l, n, \mu} \, \mcV_{l, \sigma}(\bfx) \, P_{n}^{(\gamma + k + 1 \, , l)}(\rho)  \\
\mbox{where } \ _{[x , \bfp]}A_{l, n, \mu} &= \ 
  _{\bfp}A_{l+1, n, \mu} \, (n + l + 1) \ + \ _{\bfp}A_{l-1, n+1, \mu} \, (n + \gamma + l + k + 1) \, , \ \mbox{ and }  \\
\frac{\partial}{\partial y} \partial_{\bfp} u(\bfx) &= \ 
\sum_{l, n, \mu} \, _{[y , \bfp]}A_{l, n, \mu} \, \mcV_{l, \sigma^{*}}(\bfx) \, P_{n}^{(\gamma + k + 1 \, ,  l)}(\rho)  \\
\mbox{where } \ _{[y , \bfp]}A_{l, n, \mu} &= \ \left\{ \begin{array}{l} 
(\mp) \, _{\bfp}A_{l+1, n, \mu} \, (n + l + 1) \ + \ (\pm) \, _{\bfp}A_{l-1, n+1, \mu} \, (n + \gamma + l + k + 1) \, , \ \mbox{ if } \sigma = \mu \\
(\pm) \,  _{\bfp}A_{l+1, n, \mu} \, (n + l + 1) \ + \ (\mp) \, _{\bfp}A_{l-1, n+1, \mu} \, (n + \gamma + l + k + 1) \, , \ \mbox{ if } \sigma = \mu^{*} 
\end{array} \right. .
\end{align*}

Then, using orthogonality (see \eqref{poit5p5}, \eqref{mjh0p5}, \eqref{mjh1p5}),
\begin{align*}
& Q  \:= \ \left\| \frac{\partial}{\partial x} \partial_{\bfp} u \right\|^{2}_{L^{2}_{\gamma + k + 1}} \ + \ 
    \left\| \frac{\partial}{\partial y} \partial_{\bfp} u \right\|^{2}_{L^{2}_{\gamma + k + 1}}   \\
& \ = \  \sum_{n}  \big( \, _{\bfp}A_{1 , n, 1}^{2} \ + \ _{\bfp}A_{1 , n, -1}^{2} \big) \, (n + 1)^{2}  
 \, \| \mcV_{0 , 1}(\bfx) P_{n}^{(\gamma + k + 1 \, , \, 0)}(\rho) \|^{2}_{\gamma + k + 1}    \\
 & \ \ \ + \ 2 \,  \sum_{l \ge 1, \, n, \mu}  
 \Big(  \, _{\bfp}A_{l+1 , n, \mu}^{2} (n + l + 1)^{2} \, + \, _{\bfp}A_{l-1 , n+1, \mu}^{2} \, (n + \gamma + l + k + 1)^{2} \Big)
 \, \| \mcV_{l , 1}(\bfx) P_{n}^{(\gamma + k + 1 \, , \, l)}(\rho) \|^{2}_{\gamma + k + 1}   \, .
%
\end{align*}

From Lemma \ref{lmaeq3},
\[
 \| \mcV_{l, \mu}(\bfx) \, P_{n}^{(\gamma + k + 1 \, ,  l)}(\rho) \|^{2}_{L^{2}_{\gamma + k + 1}}  \ = \ 
 \frac{2n + \gamma + l + k + 1}{2n + \gamma + l + k + 2} \, \frac{n + \gamma + k + 1}{n + \gamma + l + k + 2} \, 
  \| \mcV_{l, \mu}(\bfx) \, P_{n}^{(\gamma + k \, ,  l)}(\rho) \|^{2}_{L^{2}_{\gamma + k}} \, .
 \]
 Hence, 
 \begin{align}
 Q &\sim \ \sum_{l \ge 0, \, n, \mu} \,  _{\bfp}A_{l+1, n, \mu}^{2} \, (n + 1) \, (n + l + 2) \, 
 \| \mcV_{l, 1}(\bfx) \, P_{n}^{(\gamma + k \, ,  l)}(\rho) \|^{2}_{\gamma + k}     \nonumber \\
 & \quad \ + \ 
 \sum_{l \ge 1, \, n, \mu} \, 
 _{\bfp}A_{l-1, n+1, \mu}^{2} \,  \, (n + 2) \, (n + l + 1)  \, 
 \| \mcV_{l, 1}(\bfx) \, P_{n}^{(\gamma + k \, , l)}(\rho) \|^{2}_{\gamma + k} \, .  \label{mkl4}
 \end{align}
 
 Using Lemmas \ref{lmaeq4} and \ref{lmaeq5},
 \begin{align}
& \| \mcV_{l, \mu}(\bfx) \, P_{n}^{(\gamma + k \, ,  l)}(\rho) \|^{2}_{L^{2}_{\gamma + k}}  \nonumber \\
 & \hspace{0.2in} = \ 
 \frac{C_{l, \mu}}{C_{l+1, \mu}} \, \frac{2n + \gamma + l + k + 2}{2n + \gamma + l + k + 1} \, \, 
 \frac{n + \gamma + l + k + 1}{n + l + 1} \, 
   \| \mcV_{l+1, \mu}(\bfx) \, P_{n}^{(\gamma + k \, , \,  l+1)}(\rho) \|^{2}_{L^{2}_{\gamma + k}}   \, , \label{mkl5} \\
&  \| \mcV_{l, \mu}(\bfx) \, P_{n}^{(\gamma + k \, ,  l)}(\rho) \|^{2}_{L^{2}_{\gamma + k}}  \nonumber \\
 &  \hspace{0.2in} = \ 
 \frac{C_{l, \mu}}{C_{l-1, \mu}} \, \frac{2n + \gamma + l + k + 2}{2n + \gamma + l + k + 1} \, \, 
 \frac{n +  1}{n + \gamma + 1} \, 
   \| \mcV_{l-1, \mu}(\bfx) \, P_{n+1}^{(\gamma + k \, , \,  l-1)}(\rho) \|^{2}_{L^{2}_{\gamma + k}}   \, . \label{mkl6} 
\end{align}

Substituting \eqref{mkl5} and \eqref{mkl6} into \eqref{mkl4} and reindexing we obtain
\begin{align*}
Q  &= \  \left\| \frac{\partial}{\partial x} \partial_{\bfp} u \right\|^{2}_{L^{2}_{\gamma + k + 1}} \ + \ 
    \left\| \frac{\partial}{\partial y} \partial_{\bfp} u \right\|^{2}_{L^{2}_{\gamma + k + 1}}  \\
& \sim \   \sum_{l \ge 0, \, n, \mu}  \,  _{\bfp}A_{l+1, n, \mu}^{2} \, (n + 1) \, (n + l + 2) \, 
 \| \mcV_{l+1, 1}(\bfx) \, P_{n}^{(\gamma + k \, ,  \, l+1)}(\rho) \|^{2}_{\gamma + k}     \nonumber \\
 & \quad \ + \ 
 \sum_{l \ge 1, \, n, \mu} \, 
 _{\bfp}A_{l-1, n+1, \mu}^{2} \,  \, (n + 2) \, (n + l + 1)  \, 
 \| \mcV_{l-1, \mu}(\bfx) \, P_{n+1}^{(\gamma + k \, , \,  l-1)}(\rho) \|^{2}_{\gamma + k}  \nonumber \\
&\sim \ \sum_{2n + l \, \ge \, 1 \, , \, \mu} \,   _{\bfp}A_{l, n, \mu}^{2} \, (n + 1) \, (n + l + 1) \, 
 \| \mcV_{l, 1}(\bfx) \, P_{n}^{(\gamma + k \, ,  \, l)}(\rho) \|^{2}_{\gamma + k}   \, .
\end{align*}

Note that
\[
| u |^{2}_{H^{k+1}_{\gamma}} \ = \ 
\sum_{\bfp \in \bfX^{k}} \left( \left\| \frac{\partial}{\partial x} \partial_{\bfp} u \right\|^{2}_{L^{2}_{\gamma + k + 1}} \ + \ 
    \left\| \frac{\partial}{\partial y} \partial_{\bfp} u \right\|^{2}_{L^{2}_{\gamma + k + 1}}  \right) \, .
\]
Then, with the induction hypothesis, \eqref{mkl2}, \eqref{mkl3} 
( with $(n + 1) \, (n + l + 1) \, 
 \| \mcV_{l, 1}(\bfx) \, P_{n}^{(\gamma + k \, ,  \, l)}(\rho) \|^{2}_{\gamma + k} $ associated with
 $\| \mcV_{l, \mu}(\bfx) \,
 P_{n}^{(\gamma + k \, , \, l)}(\rho) \|^{2}_{L^{2}_{\gamma + k}}$ ),
we obtain
\begin{align*}
| u |^{2}_{H^{k+1}_{\gamma}} &\sim \ 
 \sum_{\bfp \in \bfX^{k}} \, \sum_{l, n, \mu}  \,  _{\bfp}A_{l, n, \mu}^{2} \, (n + 1) \, (n + l + 1) \, 
 \| \mcV_{l, 1}(\bfx) \, P_{n}^{(\gamma + k \, ,  \, l)}(\rho) \|^{2}_{L^{2}_{\gamma + k}}  \\
&\sim \ \sum_{2n + l \, \ge \, k + 1 \, , \, \mu} (n + 1)^{k + 1} \, (n + l + 1)^{k + 1} \, a_{l, n, \mu}^{2} \,  
\| \mcV_{l, \mu}(\bfx) \, P_{n}^{(\gamma + k \, ,  \, l)}(\rho) \|^{2}_{L^{2}_{\gamma + k}}  \, ,
\end{align*}
where in the last step we have used Lemma \ref{lmaker3}. \\
\mbox{ } \hfill \qed



\end{document}